\documentclass[]{interact}

\usepackage{epstopdf}
\usepackage[caption=false]{subfig}

\usepackage{graphicx}

\usepackage{lineno}

\usepackage{amssymb}
\usepackage{amsmath}
\usepackage{microtype}
\usepackage{url}

\newcommand{\R}{\mathbb{R}}

\newcommand{\E}{\mathbb{E}}
\renewcommand{\H}{\mathbb{H}}
\renewcommand{\S}{\mathbb{S}}
\newcommand{\T}{\mathbb{T}}
\newcommand{\SL}{\mathrm{SL}_2(\R)}
\newcommand{\USL}{\widetilde{\SL}}

\newcommand\norm[1]{\left\lVert#1\right\rVert}

\newcommand{\X}[1]{\frac{\partial}{\partial x_{#1}}}

\newcommand{\eDot}[2]{\langle#1,#2\rangle_{\mathbb{E}}}
\newcommand{\hDot}[2]{\langle#1,#2\rangle_{\mathbb{H}}}
\renewcommand{\dot}[2]{\langle#1,#2\rangle}

\newenvironment{psmallmatrix}
  {\left(\begin{smallmatrix}}
  {\end{smallmatrix}\right)}

\theoremstyle{plain}
\newtheorem{theorem}{Theorem}[section]

\newtheorem{conjecture}[theorem]{Conjecture}

\theoremstyle{definition}

\theoremstyle{remark}

\usepackage{float}

\usepackage[]{todonotes}

\begin{document}

\title{How to see the eight Thurston geometries}


\author{
\name{Tiago Novello\textsuperscript{a}, Vin\'icius da Silva\textsuperscript{b}, Luiz Velho\textsuperscript{a}, Mikhail Belolipetsky\textsuperscript{a}}
\affil{\textsuperscript{a}IMPA, Rio de Janeiro, Brazil;~\\ \textsuperscript{b}PUC-Rio, Rio de
	Janeiro, Brazil}
}

\maketitle

\begin{abstract}
In this expository paper, we present a survey about the history of the geometrization conjecture and the background material on the classification of Thurston's eight geometries. 
We also discuss recent techniques for immersive visualization of relevant three-dimensional manifolds in the context of the Geometrization Conjecture.  
\end{abstract}


\section{Introduction} \label{sec:introduction}
An $n$-\textit{manifold} is a space locally modeled on the $n$-dimensional Euclidean space. The global structure of a manifold can be non-trivial and very different from the Euclidean space. Manifolds are important objects in mathematics and physics because they allow more complicated structures to be expressed in terms of the relatively well-understood properties of simpler spaces.

This expository paper gives an overview of the problem of visualizing $3$-manifolds, which is not as easy as visualizing the Euclidean space. Specifically, the scene objects are set in $3$-manifold spaces. Recent works~\cite{ vc-rtorb-2014, global_illumination, nilsolsl2, velho2020immersive} avoid modeling perspective views by using a generalization of the \textit{ray tracing} algorithm:
a color is given to each point and \textit{tangent direction} by tracing a ray and finding its intersections with the scene objects. Recall from physics that light travels along the rays --- the paths that locally minimize length.
Tracing a ray requires \textit{geometry}, and finding geometry can be a hard task as we will see along this text. 

We can think of $3$-manifolds as spaces representing the shape of the universe, since from our view they look like the Euclidean $3$-space. This is a three-dimensional version of the fact, for example, that the surface of the earth (a $2$-sphere) is locally similar to a plane. 
For a $3$-manifold example, consider the set of points equidistant from a fixed point in the four-dimensional Euclidean space --- the \textit{$3$-sphere}. This space plays a central role in the study of $3$-manifolds being the main actor in the Poincar\'e conjecture. 

The dimension is a hard constraint on $n$-manifolds viewing; our eyes only see up to three dimensions. 
The $2$-manifolds can be visualized \textit{extrinsically} using a three-dimensional Euclidean space to illustrate its \textit{ambient space}, and \textit{intrinsically} by embedding the oriented surface in the Euclidean space and visualizing it through classical algorithms: \textit{rasterization} or \textit{ray tracing}.

\newpage

The problem of visualizing $3$-manifolds is much harder. However, in 1998, Thurston published \textit{``How to see $3$-manifolds''}~\cite{thurston1998see}, discussing the ways to visualize a  $3$-manifold using our spatial imagination and computer aid. Many tools in $3$-manifold theory are inspired by human spatial and geometrical instincts. Thus, the human mind is trained to understand the kinds of geometry that are needed for 3-manifolds. Finding a ``geometry'' for a given $3$-manifold is related to the \textit{Thurston's geometrization conjecture}, which encapsulates the Poincar\'e conjecture. We will discuss the geometrization conjecture in more details, roughly speaking the conjecture says that each $3$-manifold admits a decomposition into nice geometric pieces each modeled by one of the eight Thurston's geometries. 
This paper presents a visualization of these geometric structures. 

Since higher dimensional manifolds can not be used to visualize $3$-manifolds, an immersive approach based on a \textit{ray tracing} algorithm can be used. Rasterization is not appropriate for this scenario because perspective projection in non-Euclidean spaces is computationally nontrivial. On the other hand, a scene embedded in a $3$-manifold can be ray traced through geodesics: given a point (eye) and a direction (from the eye to the pixel) for each one. When it reaches an object its \textit{shading} can be computed. This is the visualization approach considered in our paper. 

\vspace{0.1cm}

We now review some important projects related to the visualization of $3$-manifolds.
In 1990's the Geometry Center started a program, under the leadership of William Thurston, for visualization of non-Euclidean geometries. A tool called {\em Geomview}~\cite{geomview} for the exploration of the classical geometries was developed based on OpenGL application programming interface. The group also used the \emph{virtual reality} (VR) installations for providing the sensation of being immersed inside such three-dimensional space.  Mathenautics~\cite{mathenautics} and Alice~\cite{alice} are two of their projects.

Later, Jeff Weeks proposed real-time visualization of the classical~\cite{Weeks2002} and product~\cite{Weeks06} geometries.
Recently, Hart et al.~\cite{hart2017non} presented a sophisticated immersive exploration of such curved spaces using VR, and Weeks~\cite{weeks2020, weeks2021body} improved it with a framework allowing game development. 
He also studied the problem of mapping the user's head and hands from the physical lab to the curved space and discovered that the space \textit{holonomy} would lead to violations of the \textit{coherence} between the head and hands~\cite{weeks2021body}.

The approaches mentioned above were based on the rasterization pipeline.
This technique projects the scene surfaces on the image plane located in front of the observer. This (perspective) projection is done along the rays connecting the surface points and the observer.
However, this map cannot be directly applied in $3$-manifolds with nontrivial geometry or topology.
In the case of spaces modeled by Thurston's geometries, it is common to replicate the scene in their covering spaces to ``unroll'' the manifolds with complicated topology. For example, a scene inside a torus can be ``unrolled" by considering the translated copies of the scene inside the Euclidean space. 

The ray tracing overcomes the above difficulties operating intrinsically in the $3$-manifold. 
It launches rays from the observer towards the image plane: if the ray hits a surface point, we compute a color for the point (pixel) in the image. Thus, ray tracing operates in the opposite direction to rasterization.
Berger et al.~\cite{vc-rtorb-2014} were the first to ray trace non-Euclidean spaces. Their rendering algorithm exploited programmable compute shaders and Nvidia CUDA platform to implement ray tracing on the graphics processor. This was restricted to the Euclidean and hyperbolic spaces, where rays can be modeled by straight lines.

\newpage

Recently, a framework from \cite{da2020ray, nilsolsl2, novello2020design, velho2020immersive, 33cmb} implemented on top of Nvidia RTX was introduced for real-time immersive and interactive visualization of spaces modeled by Thurston's geometries and other non-homogeneous spaces. 
Other important attempts to visualize and explore non-Euclidean geometries are presented in \cite{coulon2020nil,coulon2020sol,coulon2020ray, kopczynski2020real}.
The aforementioned approaches adopt only a local illumination model, however, \cite{global_illumination} introduced the Riemannian illumination which allows us to synthesize ``photorealistic'' inside views of non-Euclidean spaces.

The paper is organized as follows.
Section~\ref{sec:history} gives a short historical overview of the mathematical investigation in the context of the geometrization conjecture.
Section~\ref{sec:basics} reviews the basic concepts related to $3$-manifolds and ray tracing on such spaces.
Section~\ref{sec:2D_manifolds} presents some well-known results about the topological classification and geometrization of surfaces.
Section~\ref{sec:3D_manifolds} provides an overview of the geometrization of \mbox{$3$-dimensional} manifolds. 
This section also introduces \textit{Riemannian ray tracing}, a computer graphics technique for inside visualization in Riemannian $3$-manifolds.
In Section~\ref{section:eight_geometries} we present the eight Thurston's geometries and use Riemannian ray tracing to render inside views in spaces modeled by them.
The final Section~\ref{sec-final_remarks} provides some concluding remarks.

\section{Brief History}\label{sec:history}
\subsection{Henri Poincar\'e}\label{subsection:poincare}
In 1895, Henri Poincar\'e published his \textit{Analysis situs}~\cite{poincare1895analysis}, in which he presented the foundations of topology by proposing to study spaces under continuous deformations. In this case positions are not important. The main tools for topology are introduced in this paper: manifolds, homeomorphisms, homology, and the fundamental group. He also discussed how the three-dimensional geometry was real and interesting. However, there was a confusion in this paper: Poincar\'e treated homology and homotopy as equivalent~concepts.

In 1904, Poincar\'e wrote the fifth supplement~\cite{poincare1904cinquieme} to Analysis situs, where he approached three-dimensional manifolds in much more detail. This paper clarified that homology was not equivalent to homotopy already in dimension three. The construction of the \textit{Poincar\'e dodecahedron} presented there gives an example of a $3$-manifold with trivial homology but with nontrivial homotopy. An inside view of this space can be found in Section~\ref{section:classical}. In his paper Poincar\'e proposed a question: \emph{Is it possible that a compact connected $3$-manifold with trivial homotopy can be different from the $3$-sphere?} This question later became known as the \emph{Poincar\'e conjecture}.

Poincar\'e stimulated a lot of mathematical works asking whether some particular manifold exists. Works on this question were awarded three Fields medals. In 1960, Stephen Smale proved the analogue of the Poincar\'e conjecture for $n$-manifolds with $n>4$~\cite{smale2007generalized}. In 1980, Michael Freedman proved the conjecture for $4$-manifolds~\cite{freedman1982topology}. The problem in dimension three was the harderst and remained open until 2003, when Grigori Perelman proved Thurston's geometrization conjecture and consequently the Poincar\'e conjecture as a corollary~\cite{perelman2002entropy,perelman2003ricci,Perelman_2008}.


\newpage

\subsection{William P. Thurston}\label{subsection:thurston}
Thurston's work on $3$-manifolds have a geometric inclination with roots in topology. He tried to generalize the geometrization theorem of compact surfaces. This theorem states that the geometry of any compact surface can be modeled by the Euclidean, the hyperbolic, or the spherical space. Hyperbolic geometry is the most abundant because it models all surfaces with genus greater than one. In dimension three Thurston proved that five more geometries are needed while the hyperbolic geometry still~plays~the~central~role.

In 1970s, Thurston proposed the \textit{geometrization conjecture} to which he gave solid justifications (cf. \cite{thurston1982three}). It is a three-dimensional version of the geometrization theorem which states that every $3$-manifold can be cut into pieces that are geometrizable.
In dimension three Thurston proved that the conjecture holds for a large class of $3$-manifolds, the \textit{Haken manifolds}, confirming that hyperbolic geometry plays the central role. The result is known as the \textit{hyperbolization theorem}. In 1982, Thurston received a Fields medal for his contributions to $3$-manifolds. The \textit{elliptization conjecture}, the part of geometrization which deals with the spherical manifolds, was open at that time.

\vspace{-0.3cm}

\subsection{Grigori Perelman}\label{subsection:perelman}
In 2000, the Clay Institute selected seven problems in mathematics to guide mathematicians in their research, the \textit{seven Millennium Prize Problems}~\cite{jaffe2006millennium}. The Poincar\'e conjecture was one of them. The institute offered one million dollars for the first proof of each problem. They did not know that the Poincar\'e conjecture was about to be proved by Grigori Perelman as a corollary of the proof of the geometrization conjecture.

In 2003, Perelman published three papers~\cite{perelman2002entropy, perelman2003ricci, Perelman_2008} in arXiv solving the Geometrization Conjecture. He used tools from geometry and analysis. Specifically, he used the \textit{Ricci flow}, a technique introduced by Richard Hamilton to prove the Poincar\'e conjecture. Hamilton settled the conjecture for a special case when the $3$-manifold has positive \textit{Ricci curvature}. The idea is to use Ricci flow to simplify the geometry along time. However, this procedure may create \textit{singularities} since this flow expands regions with negative Ricci curvature and contracts regions of positive Ricci curvature. Hamilton suggested the use of \textit{surgery} before the manifold collapse. The procedure gives rise to a simpler manifold, and we can evolve the flow again. Perelman was able to show that this algorithm terminates and each connected component of the resulting manifold admits one of the Thurston geometries. In other words, Perelman proved the geometrization conjecture, and consequently the Poincar\'e conjecture. Several research groups around the world have verified his proof.

We will review in more detail geometrization and the Ricci flow in Section~\ref{sec:3D_manifolds}.

\vspace{-0.3cm}

\section{Basic Concepts} \label{sec:basics}

\vspace{-0.3cm}
Several concepts are needed to relate $3$-manifolds and ray tracing. We start with some definitions on topology of manifolds, then we associate a geometry to them. 

\vspace{-0.3cm}

\subsection{Topology}\label{subsection:topology}
\textit{Topology} is the branch of mathematics that studies the shape of objects modulo continuous deformation. Informally, we can stretch, twist, crumple, and bend, but not tear or paste. The $n$-\textit{manifolds} are examples of topological spaces that are locally similar to the $n$-dimensional Euclidean space. Loops are examples of $1$-manifolds, and compact surfaces are examples of $2$-manifolds. 

The first object to capture the topology of a manifold $M$ is its \textit{fundamental group} denoted $\pi_1(M)$. It records the basic information about the shape (holes) of $M$. Introduced by Poincar\'e, the fundamental group consists of equivalence classes under continuous deformation of loops based at a given point and contained in the space.
A manifold is \textit{simply connected} if its fundamental group is trivial. 
The Poincar\'e conjecture states that each compact simply connected $3$-manifold must be homeomorphic to the $3$-sphere, i.e. to have the  $3$-sphere shape. Attempts to prove the Poincar\'e conjecture led to  discovery of many manifold constructions. 

A common manifold construction is through the quotient of ``simpler'' manifolds by special groups acting on them. This is reasonable because each manifold is uniquely covered by a simply connected manifold --- the \textit{universal covering space} (see \cite[Chapter~4]{lee2010introduction}). Informally, a manifold $\widetilde{M}$ \textit{covers} a manifold $M$ if there is a map which ``evenly covers'' a neighborhood of each point in $M$.
The covering is \textit{universal} if $\widetilde{M}$ is simply connected. For example, the two-dimensional torus is covered by the Euclidean plane.
The Poincar\'e conjecture implies that if the universal covering of a compact manifold is compact, then the covering must be the sphere.
By the above discussion, we only need to consider quotients of simply connected manifolds. 

Let $M$ be a manifold and $\Gamma$ be a discrete group acting on it. The \textit{quotient manifold theorem} (Theorem 9.16 in \cite{lee2010introduction}) states that $M/\Gamma$ is a manifold when the group $\Gamma$ acts smoothly, \textit{freely}, and \textit{properly discontinuously} on $M$. Here the action of $\Gamma$ is free if it has no fixed points and it is properly discontinuous if each point $p$ admits a neighborhood $U$ such that $U\cap g(U)= \emptyset$, for all $g\in \Gamma$ different from the identity.
For example, if the quotient manifold $\mathbb{E}^2/\Gamma$ is a compact surface, then it is the torus or the \textit{Klein bottle}~(see \cite{martelli2016introduction}). 

More examples of manifolds can be constructed from the direct product, e.g. the $n$-torus $\T^n$ is the product of the circle $\S^1$ and the $(n-1)$-torus $\T^{n-1}$. 

\subsection{Geometry}\label{subsection:geometry}
In \textit{Riemannian geometry}, manifolds have a metric which allows the introduction of \textit{geodesics}: paths that locally minimize lengths. These are the ingredients for a ray tracing algorithm on manifolds. Following the notation of do~Carmo~\cite{carmo1992riemannian}, we now present the basic definitions and examples from Riemannian geometry.

Every point of an $n$-manifold $M$ admits a neighborhood homeomorphic to an open ball in $\mathbb{R}^n$, the corresponding maps are called \textit{charts}. We need the change of charts in $M$ to be \textit{differentiable}.  
Let $\textbf{x}(x_1,\ldots,x_n)$ be a chart of a neighborhood of a point $p$.  The \textit{tangent space} $T_pM$ at $p$ is the vector space spanned by the tangent vectors $\{\X{i}(p)\}$ of the coordinate curves at $p$. 
A \textit{Riemannian metric} on a manifold $M$ is a map that assigns a positive-definite inner product $\dot{\cdot}{\cdot}$ to each tangent space, such that in coordinates $\textbf{x}(x_1,\ldots,x_n)=p$ the functions $g_{ij}(x_1,\ldots,x_n):=\dot{\X{i}}{\X{j}}_p$ are smooth. 
Expressing two vectors $u,v\in T_pM$ in terms of the associated basis, that is, $u=\sum u_i \X{i}(p)$ and $v=\sum v_i \X{i}(p)$, we obtain:
\begin{equation}\label{eq:metric}
\displaystyle\dot{u}{v}_p=\sum_{i,j=1}^n\dot{\X{i}}{\X{j}}(p)u_iv_j=\sum_{i,j=1}^ng_{ij}(p)u_iv_j.
\end{equation}
\newpage

\noindent The metric $g$ is determined by the matrix $[g_{ij}]$. The pair $(M,g)$ is a \textit{Riemannian manifold}.
For examples of Riemannian $3$-manifolds we can consider the classical Euclidean, hyperbolic and spherical spaces, as well as the non-classical: $\S^2\times\R$, $\H^2\times\R$, Nil, Sol, and $\USL$ (see Section~\ref{section:eight_geometries}). All these eight geometries are \emph{homogeneous}, that is, for any two points there is an isometry sending one to another. Only Euclidean, hyperbolic and spherical spaces are \textit{isotropic}, that is, for any two vectors in the tangent space at a point there is an isometry of the manifold sending one to another. 

Let $(N, g_N)$ and $(M, g_M)$ be Riemannian manifolds of dimension $n$ and $m$, respectively. The $(n+m)$-manifold $N\times M$ admits a Riemannian metric given by $g_n+g_m$, the \textit{product metric}. Examples of the product metrics include the geometries $\S^2\times\R$ and $\H^2\times\R$.

\textit{Lie groups} are important examples of Riemannian manifolds. A Lie group is a manifold $M$ with a group structure, where the operations $(p,q)\to p\cdot q$ and $p\to p^{-1}$ are smooth.
Thus the \textit{left multiplication} by $p\in M$, given by $L_p(q)=p\cdot q$, is a smooth map.
The classical way to define a Riemannian metric on a Lie group is by fixing an inner product $\dot{\cdot}{\cdot}_e$ in the tangent space at the identity element $e$, and then extending it by the left multiplication:
\begin{equation}\label{eq:metric_lie_group}
\dot{u}{v}_p=\dot{d(L_{p^{-1}})_p(u)}{d(L_{p^{-1}})_p(v)}_e, \,\, p\in M,\,\, u,v\in T_pM. 
\end{equation}

In the geometries considered in the geometrization conjecture, only the product geometry $\S^2\times \R$ is not realized as a left invariant metric on a Lie group.

Quotients of Riemannian manifolds by discrete groups of isometries produce new Riemannian manifolds.
Specifically, the quotient $M/\Gamma$ of a Riemannian manifold $M$ by a discrete group $\Gamma$ acting isometrically on it has the \textit{geometric structure} modeled by $M$. 
This quotient corresponds to a covering, so we often consider $M$ being simply connected. There are exactly three Riemannian surfaces modeling the geometry of all closed compact surfaces (see Section~\ref{sec:surfaces}). In dimension three the list is increased by five special examples of product and non-isotropic geometries. 
These are \textit{model geometries}: complete simply connected Riemannian manifolds such that each pair of points have isometric neighborhoods.

\vspace{0.3cm}
We now define the main concept for ray tracing in this context. A \textit{geodesic} in a Riemannian manifold $(M,g)$ is a curve $\gamma(t)=(x_1(t),\ldots,x_n(t))$ with null \textit{covariant~derivative}:
\begin{equation}\label{eq:geodesic_equation}\small
	\frac{D}{dt}\gamma'=\sum_{k=1}^{n}\left(x''_k+\sum_{i,j=1}^{n}\Gamma^k_{ij}x'_ix'_j\right)\frac{\partial}{\partial x^k}=0 \Longleftrightarrow x''_k+\sum_{i,j=1}^{n}\Gamma^k_{ij}x'_ix'_j=0, \,\, k=1,\ldots,n.
\end{equation}
This differs from the classical case by the addition of $\sum\Gamma^k_{ij}x'_ix'_j$, which includes the \textit{Christoffel symbols} $\Gamma^m_{ij}$
 of $(M, g)$. 
To linearize system \eqref{eq:geodesic_equation}, we add new variables being the first derivatives $y_k = x'_k$, obtaining thus the \textit{geodesic flow} of $(M, g)$:
\begin{equation}\label{eq:geodesic_equation_bundle}
\left\{
\begin{array}{ll}
x'_k& = y_k  \\[0.0cm]
y'_k& =\displaystyle-\sum_{i,j=1}^{n}\Gamma^k_{ij}y_iy_j, \,\,\, k=1,2,\ldots,n. 
\end{array}
\right.
\end{equation}

\newpage

Let $(M,g)$ be a Riemannian $3$-manifold.
In Section \ref{sec:shader} we use the geodesic flow of $(M,g)$ to define a ray tracing algorithm in $M$. Here we give a preliminary explanation. Let $p=(x_1,x_2,x_3)\in M$ be a point (observer) and $v=(y_1,y_2,y_3)$ be a vector (pixel) in the (image) plane inside $T_pM$. This plane should be defined in ``front of the observer". Let $\gamma$ be a curve for which $(\gamma, \gamma')$ satisfies Equation~\eqref{eq:geodesic_equation_bundle} and has $(x_1,x_2,x_3,y_1,y_2,y_3)$ as the initial condition. 
Finding the intersections between the geodesic $\gamma$ and a scene in $M$ allows us to define an RGB color for the underlying pixel. To integrate the geodesic flow of $(M,g)$ we may consider numerical techniques.

\section{Two-dimensional manifolds} \label{sec:2D_manifolds}
We present some well-known results involving topology and geometry of surfaces. We assume all surfaces being compact, connected, and oriented.  Starting with the \textit{classification theorem} in terms of the \textit{connected sum}, one can represent a surface through a polygon with an appropriate edge gluing. This polygon can be embedded in one of the three two-dimensional geometry models (Euclidean, spherical, and hyperbolic). The resulting surface has the geometry modeled by one of these geometries. 

\subsection{Classification of compact surfaces}\label{sec:surfaces}
The classical way to state the classification theorem of surfaces is by the \textit{connected sum}. Removing disks $D_1$ and $D_2$ from surfaces $S_1$ and $S_2$, one obtains their connect sum $S_1\#S_2$ by identifying the boundaries $\partial D_1$ and $\partial D_2$ through a homeomorphism. 
The theorem says that any compact orientable surface is homeomorphic to a sphere or a connected sum of tori.

The proof of the classification theorem uses a computational representation of a compact surface $S$ through an appropriate pair-wise identification of edges in a polygon:
\begin{itemize}
    \item  Take a triangulation $T$ of the surface $S$; it is a well-known result that a triangulation exists;
    \item Cutting along edges in $T$ we obtain a list of triangles embedded in the plane without intersection; the edge pairing must be remembered;
    \item We label each triangle edge with a letter according to its gluing orientation;
    \item Gluing the triangles through the pairwise edge identification without leaving the plane produces a polygon $P$. The boundary $\partial P$ of $P$ can be represented as an oriented sequence of letters; 
    \item Let $a$ and $b$ be a pair of edges in $\partial P$. If the identification of $a$ and $b$ reverses the orientation of $\partial P$ we denote $b$ by $a^{-1}$, and simply $a$ otherwise;
    \item A technical result states that cutting and gluing $P$ leads us to an equivalent \textit{irreducible} polygon $Q$ with its boundary having one of following configurations:
    \begin{itemize}
        \item $aa^{-1}$, which is a sphere;
        \item $\displaystyle \sum aba^{-1}b^{-1}$, a connected sum of tori $aba^{-1}b^{-1}$.
    \end{itemize}
\end{itemize}

Although the above procedure shows that any compact connected oriented surface is homeomorphic to a sphere or a connected sum of tori, we do not know that all these surfaces are topologically different. For this, we introduce the numerical invariant called \textit{Euler characteristic}.

\newpage

Let $S$ be a compact surface and $K=(V,E,F)$ be a \textit{cell decomposition} of $S$. Here, $V$ is a set of vertices in $S$, $E$ is a set of edges (curves) with endpoints in $V$, and $T$ are the polygons (simple regions) bounded by edges in $S$.
The Euler characteristic of $S$ is the number $\chi(S)=|V|-|E|+|F|$, where $|\cdot|$ denotes the cardinality of the set. It can be shown that $\chi(S)$ does not depend on the cell decomposition $K$ of $S$ (see~Section~8~of~\cite{massey1991basic}). We will give a justification for this fact in Section~\ref{subsection:geometrization_surfaces} using the Gauss--Bonnet theorem.

As the Euler characteristic $\chi(S)$ does not depend on the cell decomposition $K$ of the compact surface $S$, we can compute it using the special cell decomposition given by the irreducible polygon $Q$ with its edge identification.
Thus, in the sphere case $aa^{-1}$, we have $\chi(S)=2-1+1=2$.
The other case is $\sum aba^{-1}b^{-1}$ which correspond to a connected sum of $g$ tori $aba^{-1}b^{-1}$, hence the Euler characteristic of the underlying surface $S$ is given by $ \chi(S)= 1-2g + 1 = 2-2g$.
Note that (in the tori case) after the edge pairing in $Q$, we get only one vertex, one face, and $2g$ edges, which justifies the formula.
Therefore, we can write the Euler characteristic of $S$ as $\chi(S)=2-2g$, where the number $g$ is also known as the \textit{genus} of $S$.

Based on the above discussion, we can classify compact connected oriented surfaces using their Euler characteristic. Therefore, this numerical invariant serves as a~dictionary.
\renewcommand{\arraystretch}{1.5}
$$\begin{array}{|c|c|}
\hline
\mathbf{\chi(S)=2-2g}& \mbox{\textbf{Surface}}  \\
\hline
\hline
2& \S^2 \\
\hline
0& \T^2 \\
\hline
-2& \T^2 \# \T^2 \\ 
\hline
-4& \T^2 \# \T^2 \# \T^2 \\
\hline
{\vdots}&{\vdots}\\
\hline
\end{array}$$

To model the geometry of these surfaces, we embed, in a special way, the irreducible polygon in one of the two-dimensional model geometries. 


\subsection{Geometrization of compact surfaces}\label{subsection:geometrization_surfaces}
We remind the well-known \textit{geometrization} theorem of compact surfaces which states that any topological surface can be modeled using only three homogeneous geometries.

\begin{theorem}[Geometrization of surfaces]
\label{theorem:geometrization_surfaces}
Any compact surface admits a geometric structure modeled by the Euclidean, the hyperbolic, or the spherical space.
\end{theorem}

This is a geometric form of the celebrated Uniformization Theorem for Riemann surfaces. We refer to a beautiful recent book by a collective of authors \cite{Uniformization-Book} that represents the history and mathematics around uniformization.

The Euclidean space $\mathbb{E}^2$ models the geometry of the $2$-torus through the quotient of $\E^2$ by the group of translations. The non-orientable Euclidean manifold is the Klein bottle. The sphere and the projective plane are  modeled by the spherical geometry. 

\newpage

For a hyperbolic surface, consider the \textit{bitorus},  which topologically is the connect sum of two tori.  The bitorus is presented as a regular polygon $P$ with $8$ sides $aba^{-1}b^{-1}cdc^{-1}d^{-1}$ as discussed above. All vertices in $P$ are identified into a unique vertex $v$. Then, the $8$ corners of $P$ are glued together producing a topological disk. Considering $P$ with the Euclidean geometry, the angular sum around $v$ equals to $6\pi$. To avoid such a problem, let $P$ be a regular polygon in the hyperbolic plane, with an appropriate scale its angles sum is $2\pi$. The edge pairing of $P$ induces a group action $\Gamma$ in the hyperbolic plane $\H^2$ such that $\H^2/\Gamma$ is the bitorus. The group $\Gamma$ tessellates $\H^2$ by the regular $8$-gons. Analogously, all surfaces represented as polygons with more than four sides are hyperbolic, implying that hyperbolic is the most abundant geometry.

The classical \textit{Gauss--Bonnet theorem} implies that these geometric models are uniquely determined by topology of the surface. Indeed, let $S$ be a smooth compact oriented surface. The Gauss--Bonnet theorem states that $\int_S K dS=2\pi\chi(S)$, where $K$ is the Gaussian curvature of $S$.
This formula connects the Gaussian curvature of $S$ (from geometry) to the Euler characteristic of $S$ (from topology).

Consider a surface $S$ modeled by a model geometry $M$.
The Gaussian curvature $K$ of $S$ coincides with the curvature of $M$, thus $K$ is constant and equals $1$, $0$, $-1$ for the spherical, Euclidean, and hyperbolic geometries, respectively. 
Therefore, the Gauss--Bonnet formula can be rewritten as $K|S|=2\pi\chi(S)$, where $|S|$ denotes the area of $S$. This formula implies some important facts.
First, $\chi(S)$ only depends on the curvature of $M$ and the area of $S$, and hence $\chi(S)$ does not depend on the cell decomposition of $S$ used on its definition.
Also, if $S$ is modeled by the Euclidean geometry, then  $K=0$ implies $\chi(S)=0$, and thus $S$ must be a torus. If $S$ is modeled by spherical geometry ($K=1$), we have $\chi(S)>0$, which shows that $S$ must be a $2$-sphere. Finally, if $S$ is modeled by the hyperbolic geometry, then $K=-1$ implies $\chi(S)<0$.
Consequently, we have that each compact oriented surface can be modeled by a unique model geometry (Euclidean, spherical, or hyperbolic). The following table shows how the Euler characteristic $\chi(S)$ determines the geometry used to model $S$.
\renewcommand{\arraystretch}{1.5}
$$\begin{array}{|c|c|}
\hline
\mathbf{\chi(S)=2-2g} &  \mbox{\textbf{Model geometry}}  \\
\hline
\hline
\chi(S)>0 & \mbox{Spherical}\\
\hline
\chi(S)=0 &\mbox{Euclidean}\\
\hline
\chi(S)<0 &\mbox{Hyperbolic}\\ 
\hline
\end{array}$$


\section{Three-dimensional manifolds} \label{sec:3D_manifolds}
It took time to develop the modern vision of manifolds in higher dimensions. For example, a version of Theorem~\ref{theorem:geometrization_surfaces} for $3$-manifolds seemed not possible until 1982, when Thurston proposed the geometrization conjecture~\cite{thurston1982three}. It states that each $3$-manifold decomposes into pieces shaped by simple geometries. There are exactly eight model geometries in dimension 3, and these are presented in more detail in Section~\ref{section:eight_geometries}. We now proceed with reviewing the geometrization.


\subsection{Classification of compact 3-manifolds}\label{section:top_3manifold}
As for surfaces, there is a combinatorial procedure to build three-dimensional manifolds from identifications of polyhedral faces. To do so, endow a finite number of polyhedra with an appropriate pair-wise identification of their faces. Each pair of faces has the same number of edges and they are mapped homeomorphically to each other. Such gluing gives a \textit{polyhedral complex} $K$, which is a $3$-manifold if and only if its Euler characteristic is equal to zero (cf. \cite[Theorem~4.3]{fomenko2013algorithmic}). 

We now take the opposite approach. Let $M$ be a compact $3$-manifold, we represent $M$ as a polyhedron $P$ endowed with a pair-wise identification of its faces. The following algorithm mimics the surface case presented in Section~\ref{sec:surfaces}.

\begin{itemize}
    \item Let $T$ be a triangulation of the manifold $M$; endorsed by the well-known triangulation~theorem;
    \item Detaching every face identification in $T$ gives rise to a collection of tetrahedra which can be embedded in $\E^3$. Remember the pairwise face gluing;
    \item Gluing in a topological way each possible coupled tetrahedra without leaving $\E^3$ produces a polyhedron $P$. The faces in the boundary $\partial P$ are pairwise identified.
\end{itemize}
The combinatorial problem of reducing the polyhedron $P$ to a standard form, as in the surface case, remains open (see page 145 in \cite{lee2010introduction}). 
Although there is no classification of compact $3$-manifolds in the sense presented for compact surfaces, it is still possible to decompose a given manifold into simpler pieces. This decomposition is not trivial. Thurston conjectured that these pieces can be modeled by one of the eight homogeneous~geometries.

The decomposition used in the geometrization theorem (to be presented in Section~\ref{section:geometrization_manifolds}) has two stages: the prime and the tori decomposition. The first is similar to the inverse of the connected sum of surfaces. It consists of cutting the $3$-manifold $M$ along a $2$-sphere such that the resulting two $3$-manifolds are not balls. A sphere 
$\S^2\subset M$ that does not bound a $3$-ball in $M$ is called \textit{essential}, and so we are cutting along essential $2$-spheres that separate $M$.
After attaching balls to the boundary of the parts, one obtains two simpler $3$-manifolds. A \textit{prime} $3$-manifold does not admit any further decomposition. One of the first results in the topology of $3$-manifolds is Kneser's theorem from 1929 showing that after a finite number of steps a manifold decomposes into prime pieces \cite{Kneser1929}, and some 30 years later Milnor proved that the decomposition is unique up to homeomorphism \cite{milnor1962unique}.

The \textit{toral decomposition} (also known as the \textit{JSJ decomposition}, named after  William Jaco, Peter Shalen, and Klaus Johannson \cite{jaco1979seifert, johannson2006homotopy}) consists of cutting a prime $3$-manifold along certain embedded tori. The result is a $3$-manifold bounded by tori that are left as boundaries, because there is no canonical way to close such holes.

\begin{figure}[ht]
    \centering
    \includegraphics[width=0.8\columnwidth]{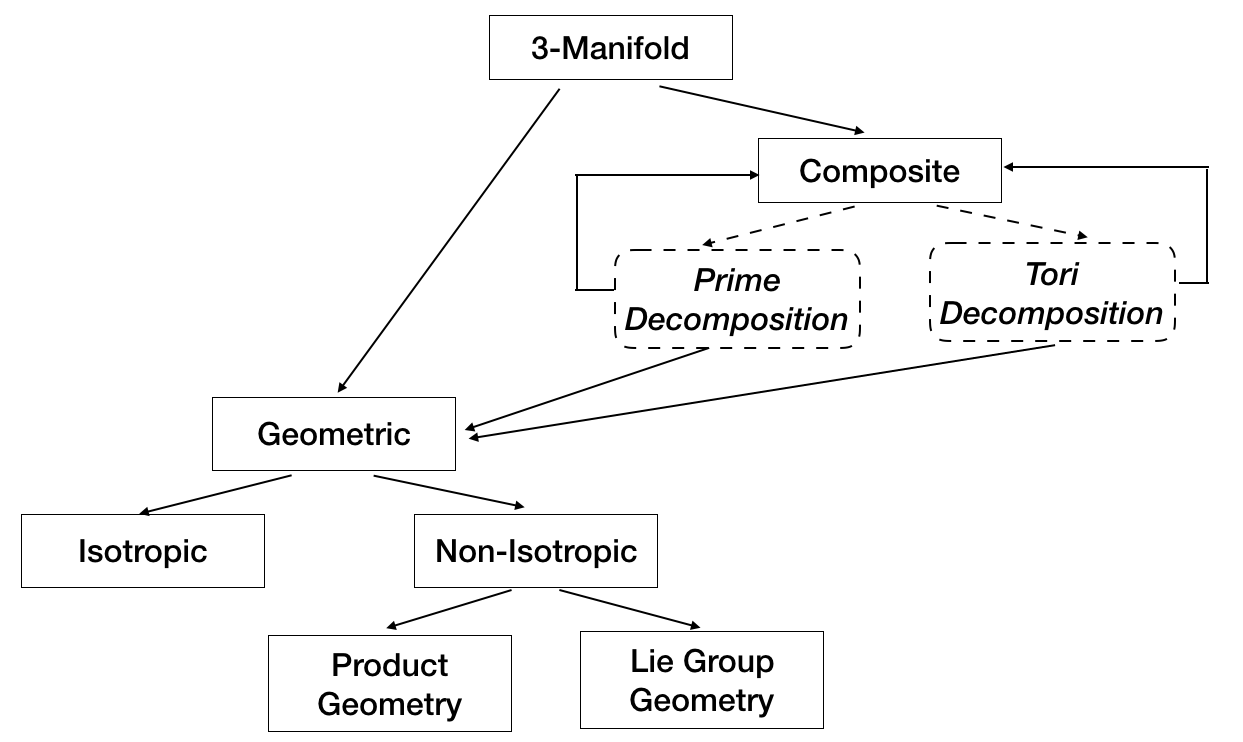}
    \caption{Evolutionary tree of a compact orientable $3$-manifold. It operates like an algorithm. The first two layers indicate the prime and tori decomposition of the $3$-manifold. The last two layers apply the geometrization~theorem.}
    \label{fig:tree}
\end{figure}

Decomposing a $3$-manifold through the above procedures produces a list of simpler manifolds, which resembles an evolutionary tree~\cite{mcmullen2011evolution}. Thurston conjectured that it is always possible to choose the separating tori so that the simplest pieces are geometric manifolds modeled by one of the eight homogeneous geometries. This is the three-dimensional case of Theorem \ref{theorem:geometrization_surfaces}--- the \textit{Thurston--Perelman Geometrization Theorem}~(see Figure~\ref{fig:tree}). 

\newpage

\subsection{Geometrization of compact 3-manifolds} \label{section:geometrization_manifolds}

To begin with we consider the three-dimensional homogeneous geometries. Now the list is longer than in dimension two but only eight Thurston's geometries admit compact quotient manifolds. These are the classical Euclidean, hyperbolic, and spherical spaces, the product spaces $\S^2\times \R$ and  $\H^2\times \R$, and the three-dimensional Lie groups \textit{Nil}, \textit{Sol}, and $\USL$. We defer the discussion of each of these eight geometries to Section~\ref{section:eight_geometries}. For now we will only use some of their basic properties. 

The quotients of the homogeneous manifolds are locally homogeneous manifolds. These are the elementary geometric pieces that we would like to see in a decomposition of a given manifold. Excluding locally homogeneous manifolds modeled on $\S^2\times \R$, all locally homogeneous $3$-manifolds have universal coverings which are topologically either $\R^3$ or $\S^3$. It is easy to check that if every embedded $2$-sphere in the universal covering of a $3$-manifold bounds a ball, then the same is true for the original manifold. Hence the only locally homogeneous $3$-manifolds with essential $2$-spheres are those whose universal cover is $\S^2 \times \R$ and there are only two such manifolds (see Section~\ref{section:eight_geometries}). Therefore, each geometric piece must be prime or contained in a prime manifold. 

The next step is to apply the toral decomposition. In 1970's Thurston conjectured that these two steps are already sufficient to get the geometric pieces:

\begin{conjecture}[Geometrization]
\label{theorem:geometrization_manifolds}\label{conj:geometrization_manifolds}
Let $M$ be a compact, orientable, prime $3$-manifold. Then there is a finite collection of disjoint, embedded, incompressible tori in $M$, so that each component of the complement admits a geometric structure modeled on one of the eight Thurston's geometries.
\end{conjecture}

Besides a good set of examples, the strongest support for the conjecture was Thurston's proof that it is true for an important class of \emph{Haken manifolds}. It is worth mentioning that Thurston's proof was influenced and somewhat resembles the proof of Andreev's theorem on classification of three-dimensional hyperbolic polyhedra~\cite{Andreev70}. We recall that an orientable, compact, irreducible $3$-manifold is called Haken if it contains an orientable, incompressible surface. Intuitively it means that the manifold is sufficiently large. The manifolds that could provide potential counterexamples to Poincar\'e's conjecture are not in this class. 

It is not hard to see that the Thurston Geometrization Conjecture implies the Poincar\'e Conjecture. Indeed, let $M$ be a compact $3$-manifold with a trivial fundamental group. If $M = M_1\# M_2\#\cdots \# M_n$ is its prime decomposition, then $\pi_1(M)$ is the free product of $\pi_1(M_i)$, therefore we have $\pi_1(M_i) = 1$ for every $i = 1, 2, \ldots, n$. Thus we can assume that $M$ is a prime manifold and the Geometrization Conjecture applies. Since $\pi_1(M) = 1$, the manifold has no incompressible tori. By Conjecture~\ref{conj:geometrization_manifolds}, $M$ must be a geometric manifold. From here the conclusion follows easily. There is only one model geometry whose space is compact and hence has compact quotient manifolds with finite fundamental group, it is the sphere $\S^3$. The condition that $M$ is simply connected implies that it is diffeomorphic to its universal cover $\S^3$. More generally, this argument shows that closed $3$-manifolds with finite fundamental groups are the three-dimensional spherical space forms (the \emph{Elliptization Conjecture}). 

We end this brief introduction to geometrization with John Morgan's~recollection~in~\cite{Morgan05}:
\begin{quote}
``It is my view that before Thurston's work on hyperbolic $3$-manifolds
and his formulation of his general Geometrization Conjecture there was no consensus among the experts as to whether the Poincar\'e Conjecture was true or false.
After Thurston’s work, notwithstanding the fact that it has no direct bearing on
the Poincar\'e Conjecture, a consensus developed that the Poincar\'e Conjecture (and
the Geometrization Conjecture) were true. Paradoxically, subsuming the Poincar\'e
Conjecture into a broader conjecture and then giving evidence, independent from
the Poincar\'e Conjecture, for the broader conjecture led to a firmer belief in the
Poincar\'e Conjecture.''
\end{quote}

\subsection{The Ricci flow}\label{section:Ricci}

It has been in the air for a long time that geometric flows can provide the key for understanding the Geometrization Conjecture. Perhaps the best intuition comes from looking at the heat flow which eventually makes the temperature to be equally distributed in a heated body. Similarly, one can expect that a certain curvature flow will make the curvature well distributed. One of the first problems encountered in this approach is which particular curvature to consider, as three dimensional manifolds have several well defined notions of curvature. Richard Hamilton was the first to investigate from this point of view the so called Ricci curvature and the corresponding Ricci curvature flow. It is not surprising that the analogy with the heat flow played important role in Hamilton's and subsequent Perelman's discoveries.  


Let $M$ be a manifold and $\{g(t)\}$ be a smooth one-parameter family of Riemannian metrics on $M$. The \emph{Ricci flow equation} is
$$\frac{d g(t)}{dt}=-2\mathrm{Ric}(g(t)),$$
where $\mathrm{Ric}(g(t))$ denotes the Ricci curvature tensor of the metric. Besides somewhat technical definition of the Ricci curvature which we will not present here, it is known to be a versatile tool with principal applications in differential geometry and general~relativity. 

The fixed points for the Ricci flow equation are the Ricci-flat manifolds. Their Ricci curvature vanishes and in dimension $3$ any such manifold is in fact flat. Next we can look at the self-similar solutions of the flow equation and we encounter the Einstein manifolds. These are the manifolds whose Ricci curvature tensor is proportional to the metric tensor. In dimension $3$ any such manifold has constant sectional curvature and thus is modeled by one of the classical geometries. These observations show a close relationship of the fixed points (up to conformal factor) of the Ricci flow in dimension $3$ and the homogeneous geometries. 

Ricci flow was introduced by Hamilton to prove that any connected closed $3$-manifold $M$ that admits a Riemannian metric with positive Ricci curvature also admits a Riemannian metric of constant positive sectional curvature. This important result establishes the Poincar\'e conjecture for the manifolds with positive Ricci curvature. 

Based on this landmark theorem, Hamilton and S.--T. Yau developed a program to attack the Poincar\'e conjecture using Ricci flow. Many important results towards this program were obtained by Hamilton and other mathematicians. A detailed description of the state of art by 1995 is given in survey \cite{Ham95}. For subsequent results we can refer to \cite{MorganTian07} and \cite{KleinerLott08}. However, after all these spectacular developments some big white spots remained wide open. Among the main mysteries were the so-called $\R\times$~\emph{cigar~soliton} singularities of Ricci flow and behavior of the flow near the points with very high~curvature. 

It is here that we encounter with the work of Grisha Perelman. In his first preprint on Ricci flow Perelman proved that the Ricci flow solution is $\kappa$-noncollapsed at sufficiently small scale and that Hamilton's cigar solitons cannot arise as a limit (thus confirming a conjecture of Hamilton). This spectacular breakthrough resolved the two big issues of Hamilton's program mentioned above. In order to prove the no local collapsing theorem Perelman introduced a new entropy functional that captures deep underlying structure properties of the Ricci flow. He then showed that local collapsing would contradict the monotonicity of the entropy.

The second Perelman's preprint presents a technically very sophisticated analysis extending the previous results to the Ricci flows with surgery, which were introduced by Hamilton in order to handle the singularities of the metrics produced by the flow. Perelman found a surgery algorithm that allowed him to control the metrics parameters and the flow after a surgery which removes a singularity is performed. This algorithm permitted Perelman to investigate the long-time behavior of the flow with surgeries and finally to prove the geometrization conjecture.  

In his work on the geometrization conjecture Perelman introduced a number of fundamental new tools which include his entropy functional, reduced length, reduced volume, and other. Since then these ingredients found new developments and new applications. One of the first applications of Perelman's work to hyperbolic geometry was given by Agol, Storm and Thurston in \cite{AgolStormThurston07}, where Perelman’s montonicity formula for the Ricci flow with surgery is used to prove a lower bound for the volume of a hyperbolic $3$-manifold with totally geodesic boundary. 

\subsection{Visualizing Riemannian $3$-manifolds}\label{sec:shader}

This section provides an introduction to \textit{Riemannian ray tracing}~\cite{global_illumination,nilsolsl2}, a computer graphics technique that allows rendering of images inside a Riemannian manifold $(M,g)$.
We will use this concept in Section~\ref{section:eight_geometries} to visualize scenes embedded in $3$-manifolds modeled by Thurston's geometries.  
We assume that light propagates along the rays in $(M,g)$ and that it is constant in vacuum.

A \textit{scene} $\mathcal{S}$ in $M$ is a set of surfaces $\{S_i\}$ inside $M$. Some of these surfaces may emit radiant energy, thus, iterating the rays leaving the light sources illuminates the scene $\mathcal{S}$. To render inside views of the illuminated scene $\mathcal{S}$ we define a camera model in $M$.

\subsubsection{Camera}
We follow the classical camera definition in~\cite[Chapter~11]{gomes2012computer} to construct a \textit{camera} in the three-dimensional manifold $M$. 
Let $p\in M$ be the camera position and $\{n,u,w\}$ be an orthogonal basis of the tangent space $T_pM$ which specifies the camera orientation. The vector $n$ is the \textit{optical} direction, $u$ is the \textit{up} direction, and $w=n\wedge u$. The reference frame $(p,\{n,u,w\})$ defines the coordinate system of the camera space. The \textit{image plane} $U$ of the camera is the plane in $T_pM$ perpendicular to the optical direction $n$ at a \textit{focal} distance $d$ from $p$. Thus, the pair $(dn,\{u,w\})$ defines a coordinate system of the image~plane $U$. 

We use the camera to render an inside view of the scene $\mathcal{S}$ in the manifold $M$. Each direction $v\in U$ can be associated with the RGB color that arrives at the observer point $p$ in the direction $-v$.
More precisely, we define a \textit{continuous image} as a map $I:U\to \mathcal{C}$ that associates an RGB color $I(v)\in\mathcal{C}$ to each direction $v\in U$. This is done by launching a geodesic ray $\gamma$, satisfying $\gamma(0)=p$ and $\gamma'(0)=v/\norm{v}$. When $\gamma$ hits a surface $S_i$ of the scene $\mathcal{S}$ at a point $q=\gamma(t_0)$, we define the RGB color $I(v)$ by computing the radiant energy emitted from $q$ in the direction $-\gamma'(t_0)$.
To compute this radiant energy, we consider the \textit{Riemannian illumination} function of $(M,g)$.

\subsubsection{Riemannian illumination}
Kajiya defined the \textit{illumination function} in the Euclidean space~\cite{kajiya1986rendering}. Novello et al. in~\cite{global_illumination} extended it to Riemannian manifolds. Specifically, let $q\in M$ be a point in a surface $S_i$ of the scene $\mathcal{S}$. The illumination function computes the amount of light emitted in a direction $v\in T_qM$. This is modeled through the integral equation over the unit hemisphere $\Omega(q)=\{v\in T_qM|\, g_q(v,N)\geq 0\mbox{ and } g_q(v,v)=1\}$, where $N$ is the normal vector of $S_i$ at $q$:
\begin{equation}\label{eq:global_Riem}
    L(q,v) =L_e(q,v) + \int_{\Omega(q)} f(q,v,w_i)L(q,w_i)g_q(w_i,N)dw_i.
\end{equation}

The \textit{bidirectional reflectance distribution function} (BRDF) $f(q,v,w_i)$ defines how the light reflects at the point $q$, and $L_e(q,v)$ is the light emitted from the surface $S_i$ at $q$ in the direction $v$. 
The term $g_q$ is the Riemannian metric $g$ restricted to the tangent~space $T_qM$. 

Considering the contribution of the directions $w_i$ leaving $q$, we divide Equation~\eqref{eq:global_Riem} into three components 
$$L(q,v) = L_e(q,v) + L_{dir}(q,v) +L_{ind}(q,v),$$ 
where $L_e(q,v)$ is the light emitted from the surface, $L_{dir}(q,v)$ is the \textit{direct contribution} coming from the light sources, and $L_{ind}(q,v)$ is the \textit{indirect contribution} reflected by other surfaces.

For computing the direct illumination, we evaluate the BRDF $f$ using unit tangent vectors. Then we perform the next event estimation, i.e. find a direction towards a light source with appropriate importance sampling. Finally, we compute the form factors, i.e. determine the size of the projected solid angle for a given subset of the manifold.

\newpage

For the indirect illumination $L_{ind}(p,v)$, we use the \textit{Monte Carlo integration} to estimate an approximation of the light reflected from other surfaces. Let $\{w_i^1,\ldots,w_i^k\}$ be $k$ vectors on the hemisphere $\Omega(q)$ chosen using a distribution density $d_w$. The Monte Carlo integration states that
\begin{equation}\label{eq:global_Riem_Approx}
    L_{ind}(q,v) \approx \frac{1}{k}\sum_{k=1}^n \frac{f(q,v,w_i^k)L(q,w_i^k)g_q(w_i^k,N)}{d_w(w_i^k)}.
\end{equation}

If $L(q,w_i^k)$ is known, the \textit{law of large number} ensures the convergence when $n\to\infty$. Otherwise, by computing these approximations for each $L(q,w_i^k)$ and repeating this procedure for a finite number of iterations, we get an approximation for the indirect illumination function $L_{ind}(q,v)$. 

\subsubsection{Image}
Let $(p,\{n,u,w\})$ be the reference frame of our camera model. We use the Riemannian illumination to define the continuous image $I:U\to \mathcal{C}$ that shows an inside view of the $3$-manifold $M$. 
In this work, we consider only an approximation of the Riemannian illumination using the direct contributions coming from the light sources, i.e. $L(q,v) \approx L_e(q,v) + L_{dir}(q,v)$. See~\cite{global_illumination} for visualizations using indirect illumination.

To represent the image $I:U\to \mathcal{C}$ in a computer, we need a discretization of its domain $U$. We use the classical computer graphics notion of image (see \cite[Chapter 6]{gomes2012computer}). 
Remember that the image plane $U$ is a plane inside $T_pM$ perpendicular to the optical direction $n$, at a focal distance $d$ from the camera position $p$. 
The pair $(dn,\{u,w\})$ defines a reference frame for the image plane $U$, so a point $dn+xu+yw\in U$ is represented as $(x,y)$ in this frame.
We use this coordinate system to discretize $U$.
First, restrict $U$ to the rectangle $[a,b]\times[c,d]$ centered at $dn$, where $a,b,c,d\in \R$. 
To discretize the rectangle consider a regular two-dimensional $k\times l$ grid $P_\Delta = \{p_{ij}=(x_i,y_j)\in\R^2\}$, where $x_i = a + i \Delta_x$, and $y_i = c + j \Delta_y$ with $\Delta_x=(b-a)/k$ and $\Delta_y=(c-d)/l$.
The grid is formed by a set of $kl$ cells $c_{ij}=[x_{i}, x_{i+1}]\times[y_{j}, y_{j+1}]$ which are called \textit{pixels} (abbreviation of \textit{picture elements}).

We sample the image function $I$ at each pixel $c_{ij}$ by considering the color $I(p_{ij})$.
This is calculated by casting a geodesic ray from the camera position $p$ towards the direction $p_{ij}$, when it hits a surface point of the scene, we compute an approximation of the Riemannian illumination at this point.
Therefore, the image function $I$ can be represented as a $k\times l$ matrix $I(p_{ij})$ which is the proper representation of an image with $k\times l$ pixels.

\vspace{0.2cm}

\begin{figure}[ht]
    \centering
        \includegraphics[width=0.486\columnwidth]{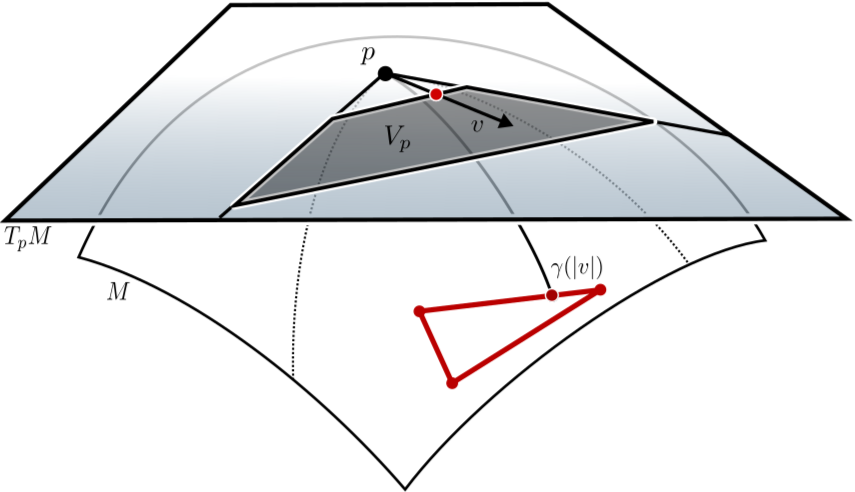}
        \includegraphics[width=0.486\columnwidth]{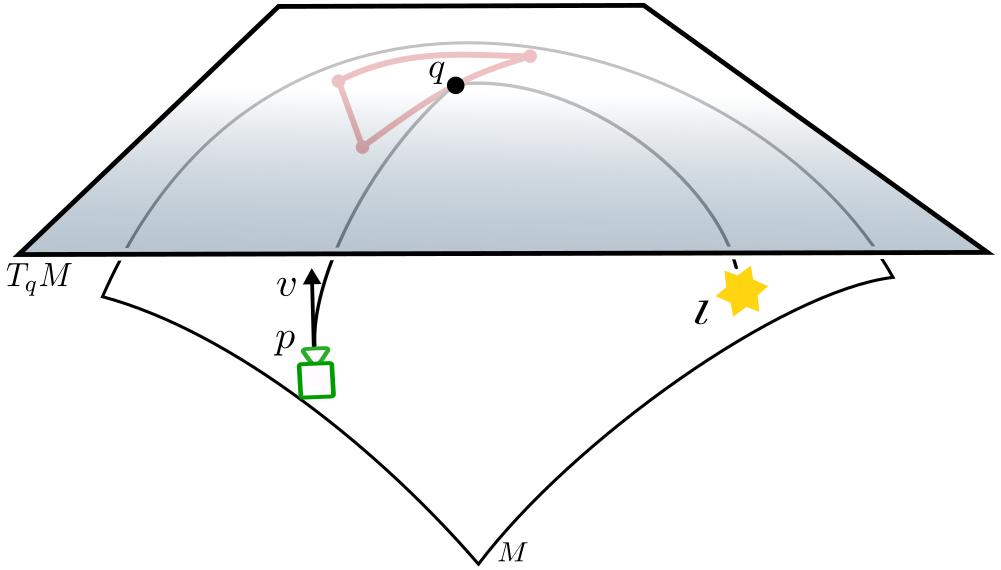}
    \caption{On the left, tracing rays in a Riemannian manifold $M$. Let $p$ be the observer and $V_p$ be the view frustum (gray region) in the tangent space $T_pM$. We launch a ray $\gamma$ towards each vector $v\in V_p$. If $\gamma$ hits a visible object (red triangle) in $\gamma(|v|)=q$ we define an RGB color for the corresponding point in the near plane of $V_p$ by considering the direct and indirect illumination (image from~\cite{nilsolsl2}). On the right, the hit point $q$ is connected with the light source $l$ through a geodesic.}
    \label{fig:frustrum_illumination}
\end{figure}

We now summarize the above visualization approach. Inside views of a scene $\mathcal{S}$ in the $3$-manifold $M$ can be rendered by tracing rays. Let $p$ be a camera position, and $P_\Delta$ be a discretization of the image plane $U$. 
We give an RGB color to each pixel in the image plane, by tracing a ray towards the associated direction.
This color is computed using an approximation of the \textit{Riemannian illumination}.
Figure~\ref{fig:frustrum_illumination} gives a schematic view of this procedure in dimension 2.

In the classical approaches ray tracing approximates physical illumination~\cite{whitted1979improved}. The above ray tracing model for Riemannian $3$-manifolds can be also used to compute a Riemannian shading function for the local~\cite{velho2020immersive} or global illumination~\cite{global_illumination}.

Visualizations of the classical Thurston geometries using simple shading for local illumination can be found in~\cite{illustratingmath}. In Section~\ref{section:classical} we describe the process in more detail and present examples. In Section~\ref{section:product}, an RGB pseudo-color based on properties of the space or attributes of the objects, such as surface normal, is used to define the Riemannian shading. The final 
Section~\ref{sec-final_remarks} presents some non-Euclidean visualizations considering global illumination (cf.~\cite{global_illumination}).

\section{The eight Thurston Geometries}\label{section:eight_geometries}

The geometrization of compact orientable $3$-manifolds provides tools for decomposing them into pieces shaped by eight homogeneous geometries. In  this section we provide the definitions and discuss some features of these geometries. We will also see why the hyperbolic geometry is the richest among the Thurston's geometries.
Excellent presentation of the eight Thurston's geometries can be found in~\cite{martelli2016introduction}, \cite{scott1983geometries}, \cite{thurston1982three}, \cite{weeks2002shape}, and other sources.

The classification of geometric pieces mentioned above makes use of the concept of \textit{Seifert manifolds}. These are closed manifolds admitting a decomposition in terms of disjoint circles. The principal two results are the following (cf. \cite{martelli2016introduction}): The first result states that a closed orientable $3$-manifold can be geometrically modeled by one of the following six geometries
$$
\R^3, \,\, \S^3, \,\, \S^2\times\R, \,\, \H^2\times\R, \,\, Nil, \text{ and } \USL
$$
if and only if it belongs to a special class of Seifert manifolds. It has a $Sol$ geometric structure if and only if it admits a particular torus bundle, called the \emph{torus \mbox{$($semi-$)$}bundle of Anosov type}. Thus the only geometry which does not have a fiber structure is the hyperbolic geometry. An outstanding conjecture of Thurston beyond the geometrization conjecture stated that hyperbolic $3$-manifolds virtually fiber over a circle. This conjecture was proved by Ian Agol in a spectacular work \cite{Agol13}. Agol's argument is based on \emph{geometric group theory}, in order to keep inline with our main topic we will refrain from discussing the details of this work.

The second result states that if a closed orientable $3$-manifold admits a geometric structure modeled by one of the eight Thurston's geometries, then it is specified by the manifold fundamental~group:
\renewcommand{\arraystretch}{1.5}
\begin{center}
\begin{tabular}{l||l}
\hline
\textbf{Fundamental group}     &  \textbf{Model geometry}\\
\hline
\hline
  Finite   & $\S^3$\\
  \hline
  Virtually cyclic & $\S^2\times \R$\\
  \hline
  Virtually abelian & $\R^3$\\
  \hline
  Virtually nilpotent & $Nil$\\
  \hline
  Virtually solvable & $Sol$\\
  \hline
  Contains a normal cyclic subgroup&\hspace{-0.25cm}
  \begin{tabular}{l|l}
      Quotient lifts\\ a finite-index subgroup & $\H^2\times \R$  \\
      \hline
      Otherwise & $\USL$
  \end{tabular}\\
  \hline
 Otherwise & $\H^3$\\
 \hline
\end{tabular}
\end{center}

The seven classes of fundamental groups mentioned above represent a restricted portion of the set of all possible fundamental groups. This implies that the hyperbolic manifolds are more abundant. We skip the precise algebraic definition of the groups involved in the table, any text on group theory can serve as an introduction.

Thurston's geometries can be divided in two classes. The isotropic geometries (Euclidean, spherical, and hyperbolic spaces) are called \textit{classical}. The \textit{product} geometries $\S^2\times \R$ and $\H^2\times \R$, and also  $Nil$, $Sol$ and $\USL$ are the \textit{non-isotropic} geometries. All these geometries are homogeneous. The classical geometries admit constant sectional curvature since they are isotropic~\cite{carmo1992riemannian}.

\subsection{Isotropic geometries}\label{section:classical}
For any dimension $n\geq 2$ there exists a unique complete, simply connected Riemannian manifold having constant sectional curvature $1$, $0$, or $-1$. These are the isotropic geometries: the sphere, the Euclidean space, and the hyperbolic space. Conversely, if a complete manifold has constant sectional curvature $1$, $0$, or $-1$, then it must be a quotient of such model geometries by a discrete group (see Proposition 4.3 in~\cite{carmo1992riemannian}).

This subsection presents the isotropic geometries, examples of manifolds modeled by them, and inside views of these spaces using non-Euclidean ray tracing techniques~\cite{velho2020immersive}. The visualizations consider only the scenes composed of Lambertian surfaces, i.e. the radiance is independent of the angle from which the surface is viewed and of the angle from which it is illuminated. The scenes are ray traced using simple shadings and no indirect light transport is considered. 
Edges of the fundamental domain of the spaces are added to highlight their structure.
For ``photorealistic'' views considering more complex scenes and indirect illumination see~\cite{global_illumination}.
For visualizations using classical rasterization techniques we refer to Weeks~\cite{Weeks2002}.

\subsubsection{Euclidean space}
\begin{figure}[H]
    \centering
    \begin{tabular}{cc@{\hskip 0.3in}c}
        \includegraphics[width=0.12\columnwidth]{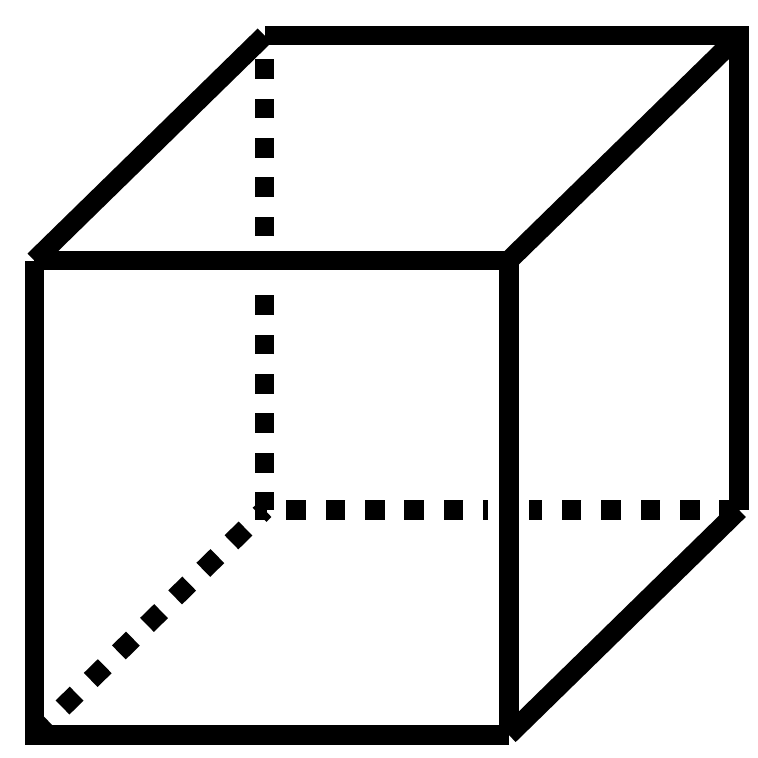}&
        \includegraphics[width=0.12\columnwidth]{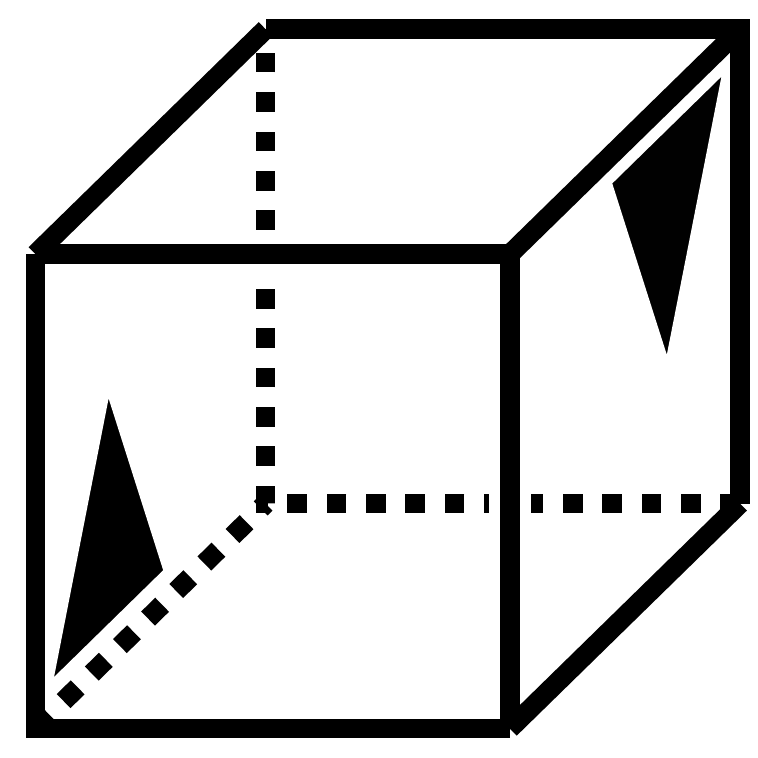}&
        \includegraphics[width=0.12\columnwidth]{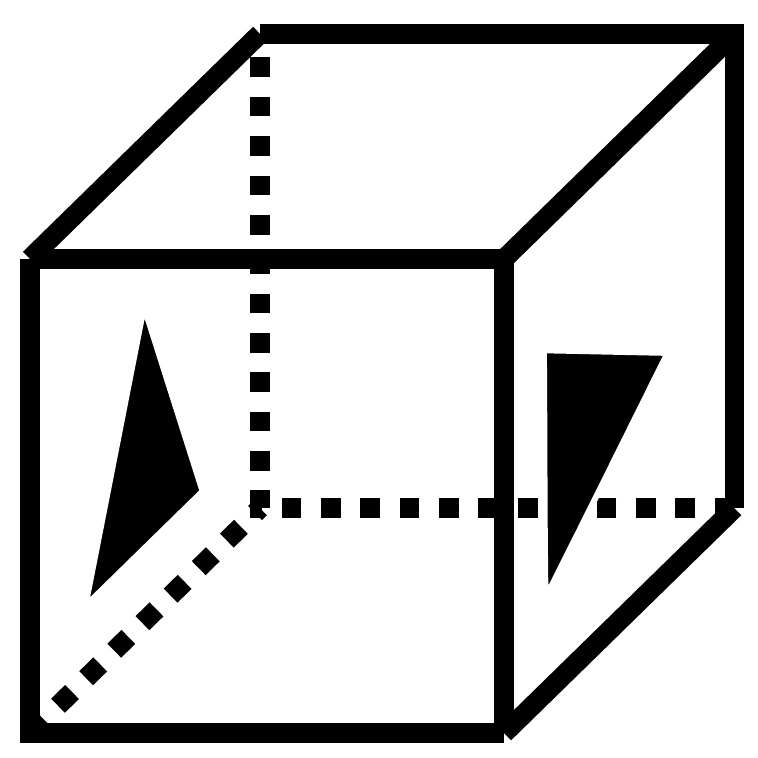}\\
        \raisebox{.4\height}{\includegraphics[width=0.14\columnwidth]{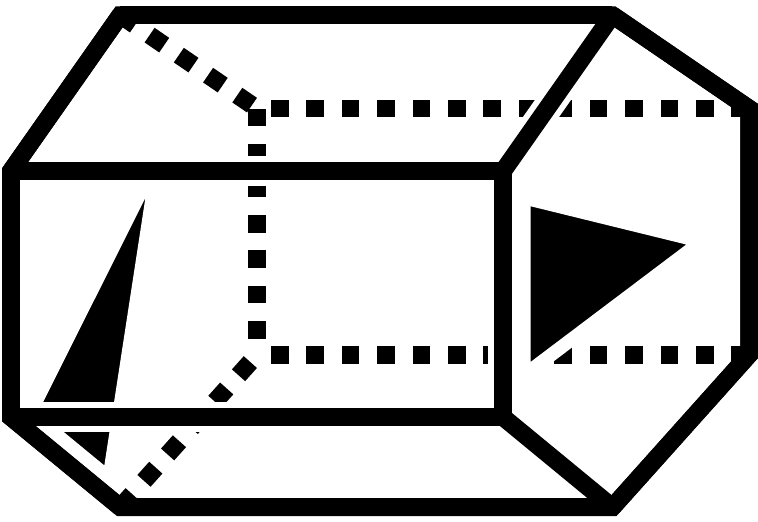}}&
        \raisebox{.4\height}{\includegraphics[width=0.14\columnwidth]{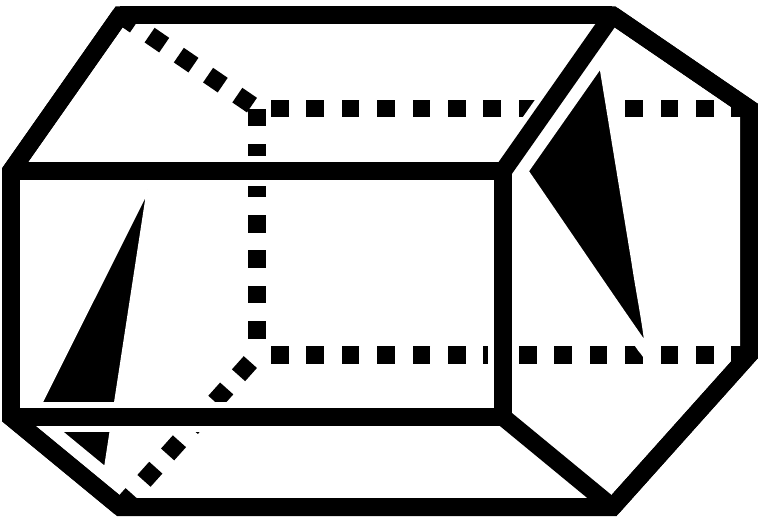}}&
        \includegraphics[width=0.12\columnwidth]{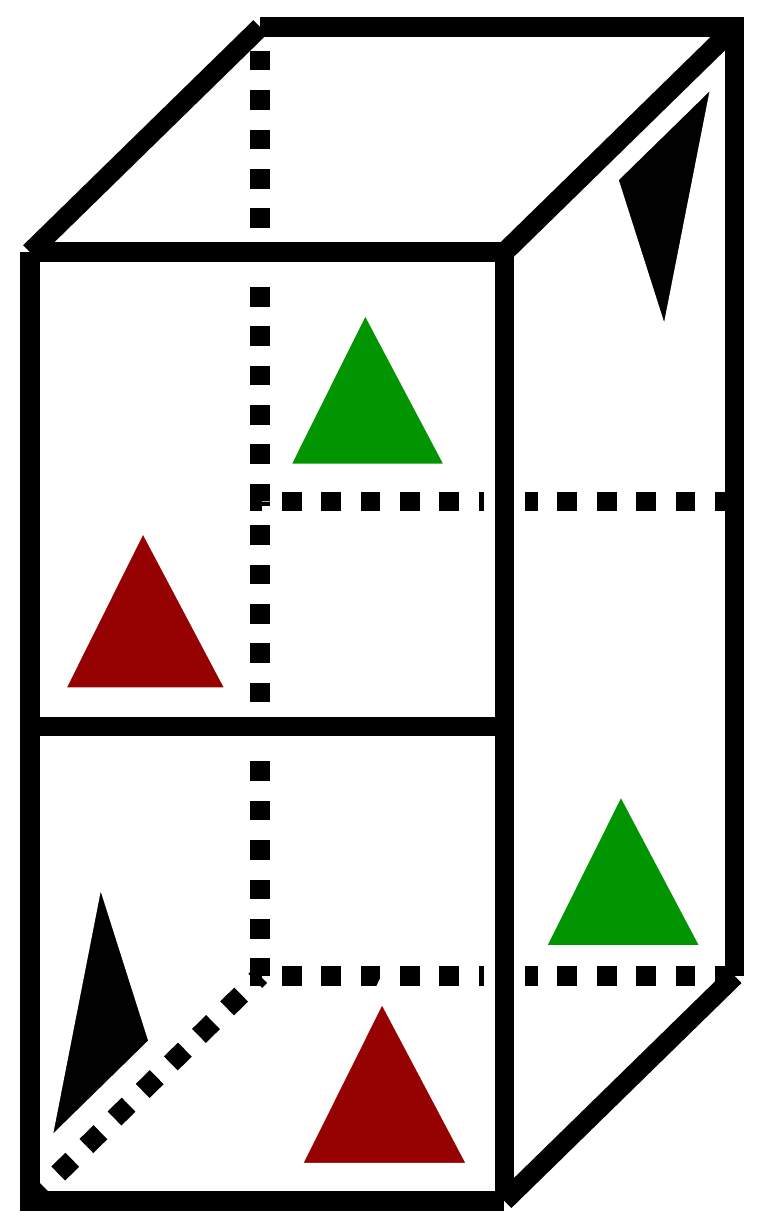}
    \end{tabular}
    \caption{The six compact oriented flat manifolds. These are built through pair-wise gluing: the faces are identified isometrically according to their labels or glued to the opposite in an obvious way.}
    \label{fig:six_flat_manifolds}
\end{figure}

In dimension two, every orientation preserving isometry of the Euclidean space $\E^2$ which has no fixed points is a translation. Therefore, if $\E^2/\Gamma$ is a compact orientable surface, it must be the torus (see \cite[Section~6.2]{martelli2016introduction}). 
In dimension three this list is increased by five more orientable manifolds since we can compose rotations with translations. These manifolds are schematically presented in Figure~\ref{fig:six_flat_manifolds}.

The \textit{Euclidean space} $\mathbb{E}^3$ is $\R^3$ with the \textit{inner product} $\langle u,v \rangle_{\E}=u_xv_x+~u_yv_y+~u_zv_z$, where $u$ and $v$ are vectors in $\R^3$. 
The \textit{distance} between two points $p$ and $q$ is $d_\mathbb{E}(p,q)~=~\sqrt{\langle p-q,p-q \rangle_{\E}}$. 
The curve $\gamma(t)=p+t\cdot v$ describes a \textit{ray} leaving a point $p$ in a direction $v$.
Analogously, for any natural $n$ we can construct the $n$-dimensional Euclidean space $\mathbb{E}^n$.

For an example of a $3$-manifold modeled by $\E^3$, consider the \textit{flat} \textit{torus} $\T^3$, obtained by gluing opposite faces of the unit cube in $\E^3$. 
The torus $\mathbb{T}^3$ is also the quotient of $\E^3$ by its group of translations generated by $(x,y,z)\to(x + 1,y,z)$, $(x,y,z)\to(x,y + 1,z)$ and $(x,y,z)\to(x,y,z + 1)$. The unit cube is the fundamental domain. 

A ray leaving a point $p\in \mathbb{T}^3$ in a direction $v$ is parameterized as $r(t)=p+t\cdot v$. For each intersection between $r$ and a face $F$ of the unit cube we update the point $p$ by $p-n$ in the opposite face, where $n$ is the normal to $F$. The direction $v$ does not need to be updated.

\begin{figure}[ht]
\centering
\includegraphics[width=0.8\columnwidth]{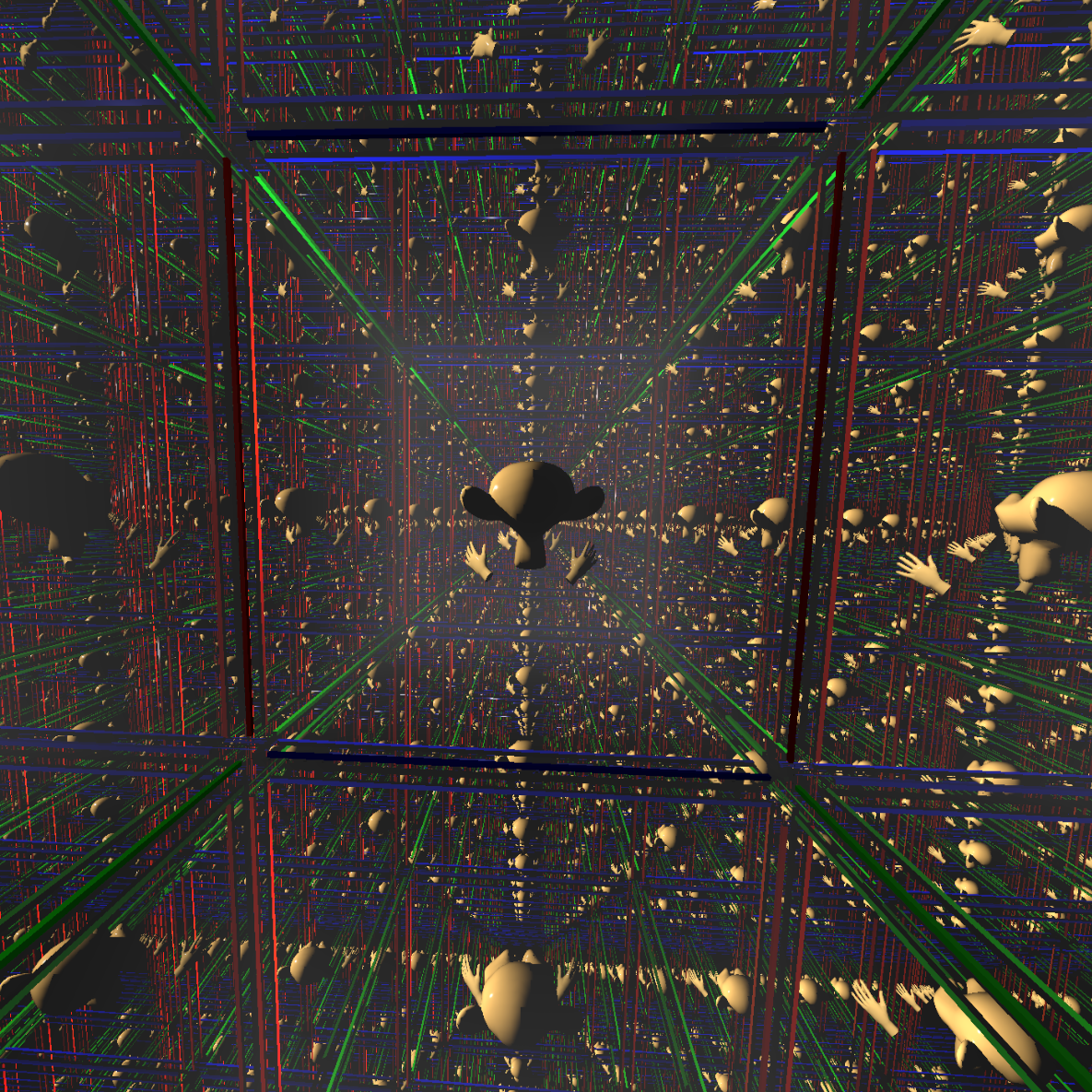}
\caption{Inside view of the flat torus. We use the cube to set up our scene: a unique mesh (Suzanne) with hands and the cube's edges. The face pairing makes the rays that leave a face to return from its opposite.}
\label{f-3torus}
\end{figure}

Thus, we have everything what is needed for an inside view of the three-dimensional torus $\T^3$. The fundamental domain receives the scene and the rays in $\T^3$ can return to it, resulting in many copies of the scene. The immersive perception is $\E^3$ tessellated by the scene copies; see Figure~\ref{f-3torus}.

Besides the torus there are exactly five more compact oriented $3$-manifold with geometry modeled by the Euclidean space, see Figure~\ref{fig:six_flat_manifolds}.

\begin{figure}[H]
\centering
\includegraphics[width=0.475\columnwidth]{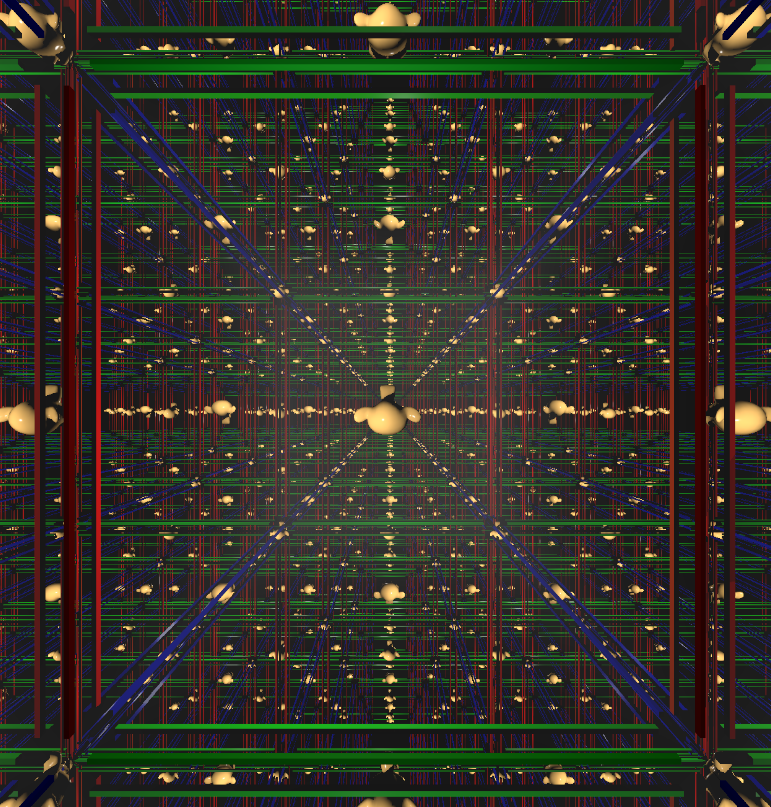}
\includegraphics[width=0.47\columnwidth]{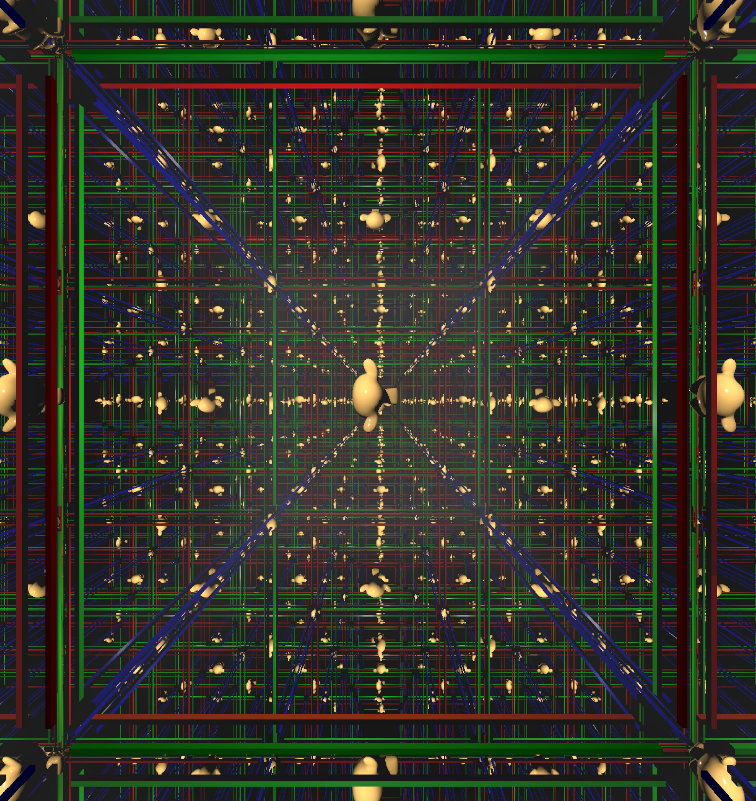}
\caption{Inside view of the flat $3$-manifolds provided by the last two examples in the first row of Figure~\ref{fig:six_flat_manifolds}.}
\label{f-flat_3_manifolds_1}
\end{figure}

The first line of examples in Figure~\ref{fig:six_flat_manifolds} presents three flat $3$-manifolds having the cube as their fundamental domains. The first example is the flat torus (illustrated in Figure~\ref{f-3torus}) and the other two spaces are similar to the torus but with a different face identification.
Figure~\ref{f-flat_3_manifolds_1} shows the immersive visualizations of these two $3$-manifolds.
On the left, one of the face identifications has a rotation of $\pi$, therefore, when looking towards that face, we see ourselves upside down in the next copy. Suzanne model was attached to the camera position to give the user a better perception of the space.
On the right, the same face identification has a rotation of $\pi/2$.

\begin{figure}[H]
\centering
\includegraphics[width=0.397\columnwidth]{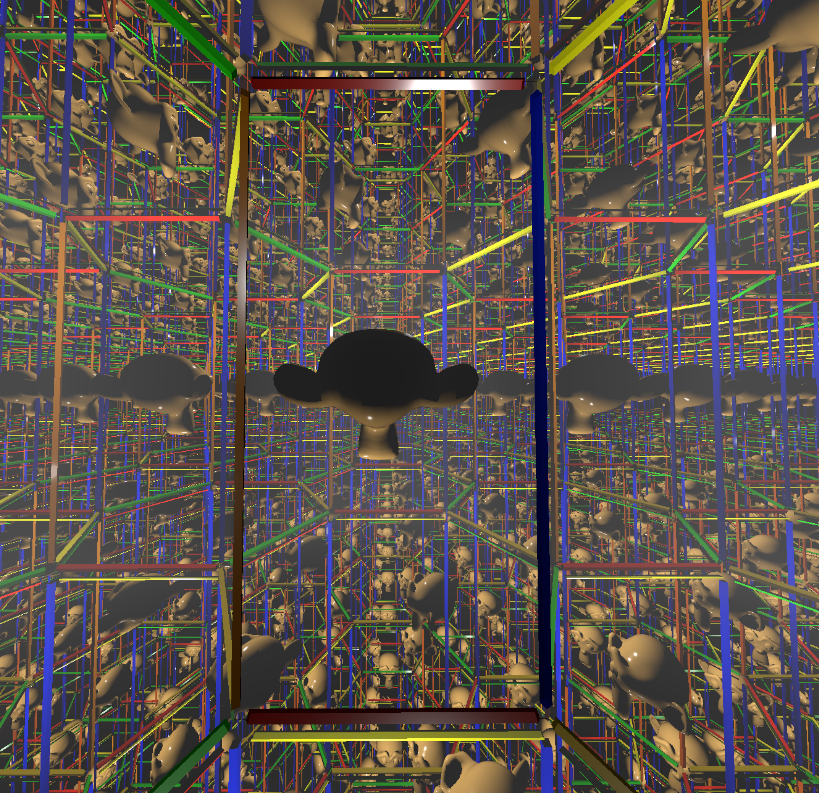}
\includegraphics[width=0.379\columnwidth]{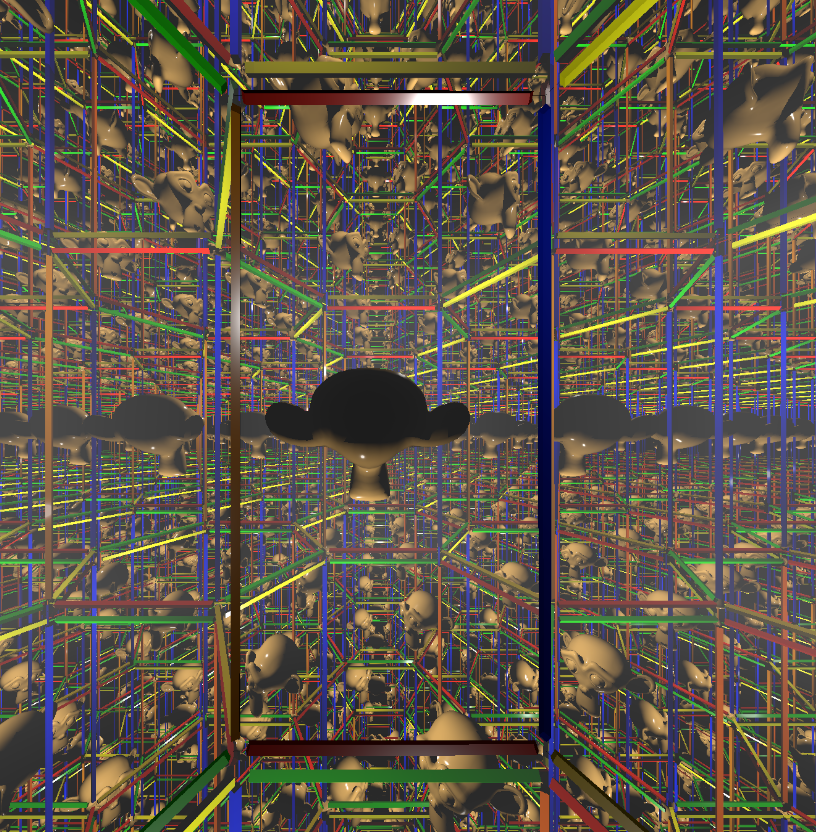}\vspace{0.05cm}
\includegraphics[width=0.79\columnwidth]{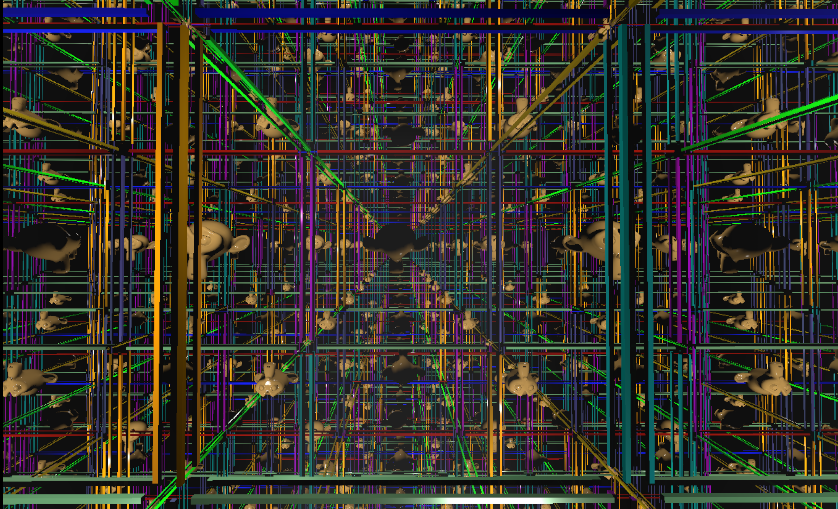}
\caption{Inside view of the flat $3$-manifolds provided by the three examples in the second row of Figure~\ref{fig:six_flat_manifolds}.}
\label{f-flat_3_manifolds_2}
\end{figure}

Figure~\ref{f-flat_3_manifolds_2} provides the inside views of the other three flat $3$-manifolds given in Figure~\ref{fig:six_flat_manifolds}. 
The first two manifolds have the hexagonal prism as their fundamental domains. The only difference in these spaces is the identification of the hexagonal faces. 
The last image gives the immersive view of the $3$-manifold modeled by the parallelepiped in Figure~\ref{fig:six_flat_manifolds}. We are looking towards the face with the green triangle label.

\subsubsection{Hyperbolic space}
Here we describe the \textit{hyperboloid} and \textit{Klein models} and present a manifold modeled by the hyperbolic $3$-space. There are plenty of hyperbolic manifolds, making this concept a central actor in the topology of $3$-manifolds (cf.~\cite{martelli2016introduction}).

The \textit{Lorentzian space} is $\R^4$ with the product $\hDot{u}{v}=u_x v_x+u_y v_y+u_z v_z-u_w v_w$, where $\{v,u\}\subset\R^4$.
The \textit{hyperbolic space} $\mathbb{H}^3$ is the hyperboloid $\{p\in \R^4|\,\,\hDot{p}{p}=-1\}$ endowed with the metric $d_{\H}(p,q)=\cosh^{-1}(-\hDot{p}{q})$, where $p$ and $q$ are points in $\mathbb{H}^3$. Due to its similarity to the sphere definition, $\H^3$ is known as the \textit{pseudo-sphere}.

A tangent vector $v$ to a point $p$ in $\mathbb{H}^3$ satisfies $\hDot{p}{v}=0$. Moreover, the \textit{tangent space} $T_p\H^3$ coincides with the set $\{v\in\R^4|\,\, \hDot{p}{v}=0\}$. The Lorentzian inner product is positive on each tangent space.

\textit{Rays} in $\mathbb{H}^3$ arise from intersections between $\mathbb{H}^3$ and planes in $\R^4$ containing the~origin. A ray leaving a point $p\in\mathbb{H}^3$ in a tangent direction $v$ is the intersection~between~$\mathbb{H}^3$ and the plane spanned by the vectors $v$ and $p$. Its parameterization is $r(t)=\cosh(t)p+\sinh(t)v$. Thus, rays in this model of the hyperbolic space can not be straight lines.

It is possible to model the hyperbolic $3$-space in the unit open ball in $\R^3$ --- known as the \textit{Klein model} $\mathbb{K}^3$--- so that in this model the rays are straight lines. More precisely, each point $p\in\H^3$ is projected in the space $\{(x,y,z,w)\in\R^4| \,\, w=1\}$ by considering $p/p_w$, the space $\mathbb{K}^3$ is obtained by forgetting the coordinate $w$. 

The hyperbolic space is an example of a \textit{non-Euclidean} geometry, since it does not satisfy the Parallel Postulate: given a ray $r$ and a point $p\notin r$, in the plane though $p$ and $r$ there exists a unique line parallel to $r$. For a ray $r$ in the hyperbolic space and a point $p\notin r$ there is an infinite number of lines in the plane $(p,r)$ which go through $p$ and do not intersect $r$.

For a compact $3$-manifold modeled by hyperbolic geometry consider the \textit{Seifert--Weber dodecahedral space}. It is the dodecahedron with an identification of its opposite faces with a clockwise rotation of $3\pi/10$. 
The face pairing groups edges into six groups of five, making it impossible to use Euclidean geometry: the regular Euclidean dodecahedron has a dihedral angle of about $116$ degrees. The desired dodecahedron should have a dihedral angle of $72$ degrees, which is possible in the hyperbolic space by considering a
regular dodecahedron of an appropriate diameter.

\begin{figure}[ht]
\centering
\includegraphics[width=0.72\columnwidth]{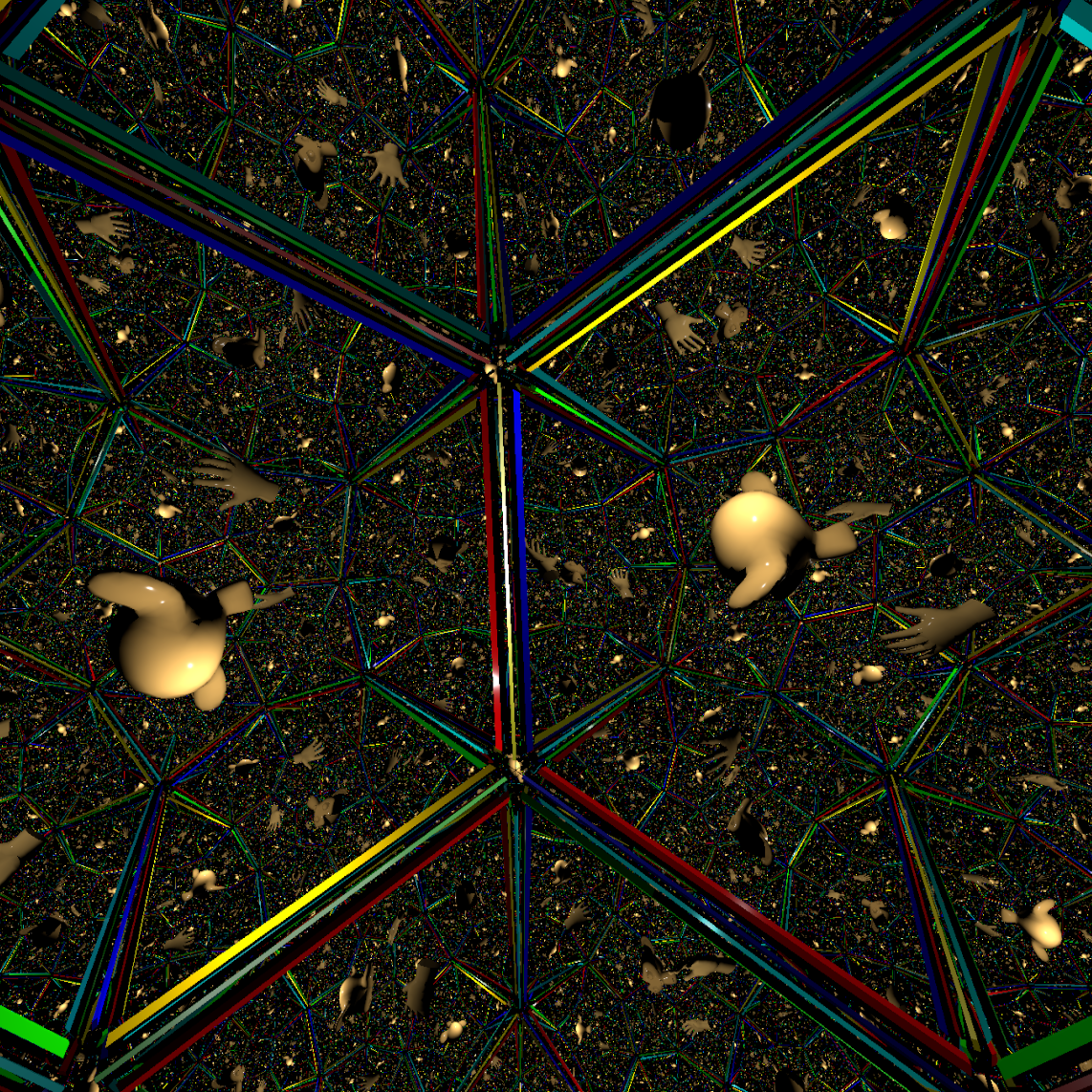}
\caption{Immersive view of the Seifert--Weber dodecahedron. We use the dodecahedron to set up our scene: a unique Suzanne with hands and dodecahedron's edges. The face pairing make the rays that leave a face to return, with an additional rotation, from its opposite face, giving rise to many copies of the scene.}
\label{f-sw-dodeca}
\end{figure}

Now we can ray trace the Seifert--Weber dodecahedron. A ray leaving a point $p\in M$ in a direction $v$ is given by $r(t)=p+tv$ since we are using the Klein model. For each intersection between $r$ and a dodecahedron face, we update $p$ and $v$ through the hyperbolic isometry that produces the face pairing above. This isometry is quite distinct from Euclidean isometries~\cite{Gunn93}.
The immersive perception of $M$ using this approach is a tessellation of $\H^3$ by dodecahedra with a dihedral angle of $ 72 $ degrees.
Figure~\ref{f-sw-dodeca} illustrates an inside view of the Seifert--Weber dodecahedral space.

\subsubsection{Spherical space}
The \textit{$3$-sphere} $\S^3$ is $\{p\in \mathbb{E}^4|\,\,\eDot{p}{p}=1\}$ with the metric $d_{\mathbb{S}}(p,q)=\cos^{-1}{\langle p,q \rangle}_\E$.
As in the hyperbolic case, a vector $v$ tangent at a point $p\in\S^3$ must satisfy $\eDot{p}{v}=0$. Therefore, the \textit{tangent space} $T_p\S^3$ at $p$ corresponds to the set $\{v\in \mathbb{E}^4|\,\eDot{p}{v}=0\}$. The space $T_p\S^3$ inherits the Euclidean inner product from $\E^4$.

A \textit{ray} in $\mathbb{S}^3$ leaving a point $p$ towards a direction $v\in T_p\S^3$ is the arc given by the intersection of $\mathbb{S}^3$ and the plane $(p,v)\subset \R^4$ spanned by $p$ and $v$. Its parameterization is $r(t)=\cos(t)p+\sin(t)v$. Therefore, spherical rays are not straight lines.

The $3$-sphere $\mathbb{S}^3$ is also a model of non-Euclidean geometry since it also fails the parallel postulate: there are no parallel lines. 
In a $2$-sphere $\mathbb{S}^2\subset \mathbb{S}^3$, any two distinct rays will intersect at exactly two antipodal points.

\begin{figure}[ht]
\centering
\includegraphics[width=4.3in]{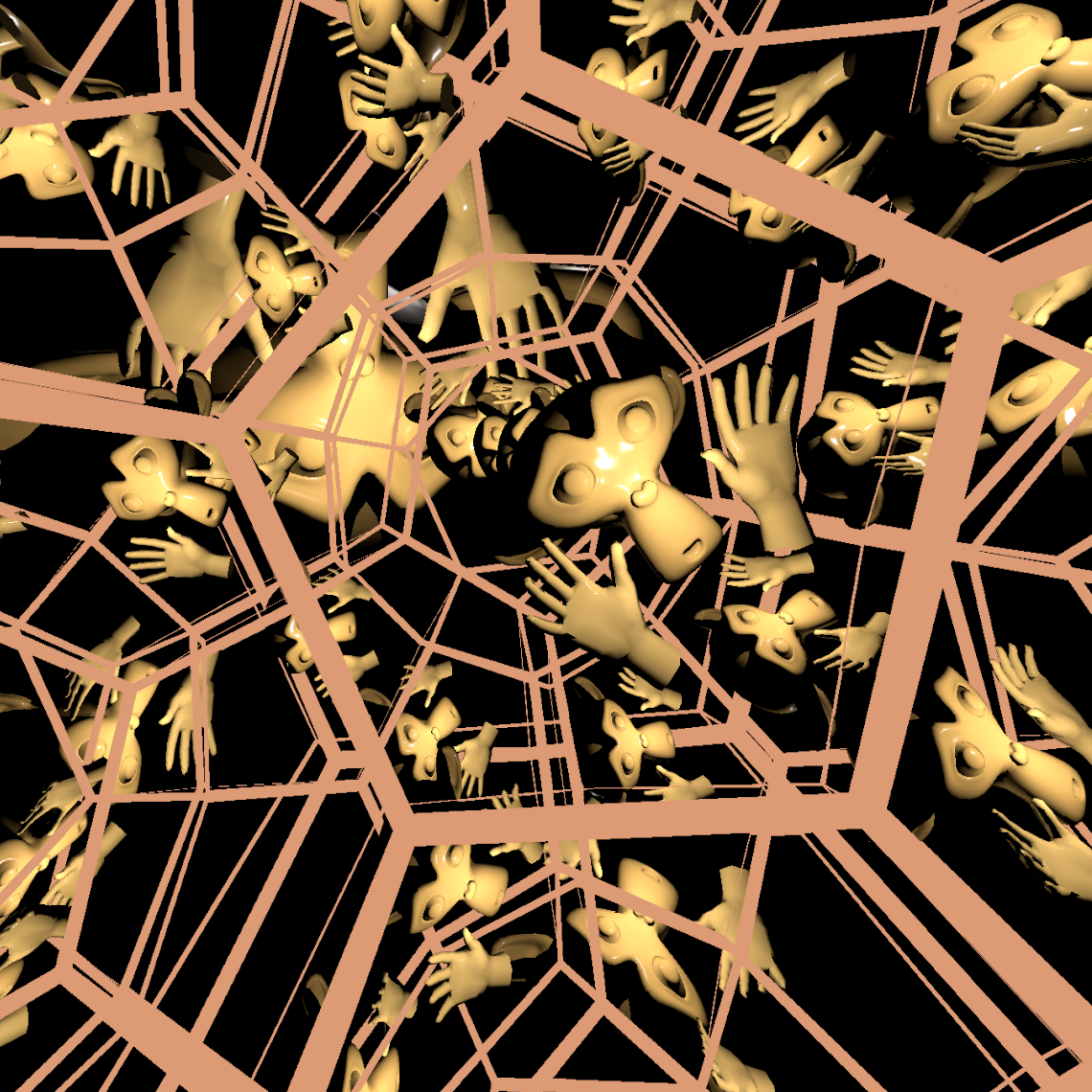}
\caption{Inside view of the Poincar\'e dodecahedron. The Suzanne with hands and the dodecahedron's edges composes the scene. The faces pairing make the rays iterate tessellating the sphere, the $120$\textit{-cell.}}
\label{f-poin-dodeca}
\end{figure}

We consider the \textit{Poincar\'e dodecahedral space} for a compact $3$-manifold modeled by spherical geometry. It is obtained by gluing the opposite faces of a dodecahedron with a clockwise rotation of $\pi/5$. The face pairing groups edges into ten groups of three, again, we cannot use Euclidean geometry. We need a dodecahedron with a dihedral angle of $120$ degrees, for this, we consider a spherical dodecahedron with an appropriate diameter. Figure~\ref{f-poin-dodeca} illustrates an inside view of the Poincar\'e dodecahedral space.

First constructed by Henri Poincar\'e, the Poincar\'e dodecahedral space (also known as the  \textit{Poincar\'e homological sphere}) has trivial homology and its fundamental group is the \textit{binary icosahedral group} of order $120$. This space is the only homology $3$-sphere, besides $\mathbb{S}^3$ itself, with a finite fundamental group.

The immersive visualization of the Poincar\'e dodecahedral space is a tessellation of $\S^3$ by $120$ dodecahedra. Each dodecahedron corresponds to an element of the binary icosahedral group. This is the boundary of the four-dimensional regular polytope, with \textit{Schl\" afli symbol} $\{5,3,3\}$, known as the $120$-cell (see
Figure~\ref{f-poin-dodeca}).

The boundary of the four-dimensional $120$-cell is a three-dimensional cellular decomposition of the $3$-sphere consisting of $120$ dodecahedra with $4$ meeting at each vertex.
The $120$-cell can be interpreted as the four-dimensional extension of the regular dodecahedron, which has Schl\" afli symbol $\{5,3\}$. The boundary of the dodecahedron is a cellular decomposition of the $2$-sphere in $12$ pentagons, with $3$ around each vertex. On the other hand, the boundary of the $120$-cell is a cellular decomposition of the $3$-sphere with $120$ dodecahedra, with $3$ around each edge.

\subsection{Non-isotropic geometries}\label{section:product}
The eight three-dimensional geometries include products of
lower-dimensional geometries, which are $\S^2\times \R$ and $\H^2\times \R$ endowed with the product metric.
We do not present immersive visualization of them because they model only few manifolds (see ~\cite[Section~12.4]{martelli2016introduction}).
Visualizations of these geometries are given by Weeks in \cite{Weeks06}.

\subsubsection{The space $\S^2\times \R$}

The geometry $\S^2\times \R$ models very few manifolds. The sectional curvature is $1$ along the horizontal directions and $0$ along the verticals. Recall that sectional curvature of a plane is the Gauss curvature associated with the surface generated by such a plane. The geometry $\S^2\times \R$ is the only Thurston geometry which is not a Lie group. It is also the only geometry that models $3$-manifolds with essential $2$-spheres.

For an example of a $3$-manifold modeled by $\S^2\times \R$, consider the $3$-manifold $\S^2\times \S$ endowed with the product metric. The geometry of $\S^2\times \S$ can not be modeled by classical geometries, since $\S^2\times \S$ has $\S^2\times \R$ as it universal covering and this is not an isotropic space. The only other example of a manifold modeled by $\S^2\times \R$ is the connected sum of projective spaces $\R P^3\#\R P^3$.

\subsubsection{The space $\H^2\times \R$}

The geometry $\H^2\times \R$ is given by the product of the hyperbolic and Euclidean metrics. Analogously to the $\S^2\times \R$ geometry horizontal and vertical planes have sectional curvature $-1$ and $0$, respectively. Compact manifolds modeled by $\H^2\times \R$ are either of the form $\Sigma\times\S$, where $\Sigma$ is a closed hyperbolic surface, or are finitely covered by such~manifolds.  

\bigskip

In \cite{Rosen2002}, Harold Rosenberg initiated systematic study of minimal surfaces in $\S^2\times\R$ and $\H^2\times \R$. This study has grown into an active area of modern differential geometry with important general theorems and beautiful examples. The product geometries provide a bridge from $\R^3$ to the other Thurston's geometries, for example $\H^3$, which came into focus in more recent investigations. Visualization of minimal and constant curvature surfaces is a feasible next step for the visualization project and we expect that here the product geometries will play a prominent role.

\newpage

The remaining three non-isotropic geometries to analyze are not products, but they admit a kind of ``bundle structure''. The first attempt to visualize these spaces in real-time (using VR) took place in 2019~\cite{nilsolsl2}. We give a brief introduction to these geometries, along with visualizations of their abstract landscapes using Riemannian ray tracing. The visualization uses RGB pseudo-color based either on the properties of the space or on the attributes of the objects, such as the surface normal, to define the Riemannian shading function. See \cite{coulon2020nil,coulon2020sol,coulon2020ray,kopczynski2020real,zenorogue, skrodzki2020illustrations} for other approaches to this~problem.

\subsubsection{Nil space}\label{sec:Nil}

Nil geometry is an $\R$-bundle over $\R^2$. This geometry is constructed from the Lie group $H$ called the Heisenberg group~\cite{martelli2016introduction}.

The \textit{Nil space} ($Nil$) is an example of a Lie group consisting of all $3\times 3$ real matrices
$$\small\left[
\begin{array}{ccc}
1&x&z\\
0&1&y\\
0&0&1
\end{array}
\right]$$
with the matrix multiplication operation. There is a natural identification of $\R^3$~with~$Nil$.

The multiplication of elements
$
(x,y,z)\cdot (x',y',z')=(x+x',y+y',z+z'+xy')
$  in $Nil$ is the sum of the coordinates with an additional term in the last one. This term makes all the difference, since in order to define a geometry in $Nil$ we consider the left multiplications $(x,y,z)\to p\cdot (x,y,z)$, for all $p\in Nil$, being isometries.

We construct a metric in $Nil$ by considering the Euclidean product in the tangent space at $e = (0,0,0)$. Then we extend it by the left multiplication.
After some calculations we obtain the scalar product between the tangent vectors $u$ and $v$ at a point $p$:
$$\small
\dot{u}{v}_p=u^T\left[
\begin{array}{ccc}
1&0&0\\
0&p_x^2+1&-p_x\\
0&-p_x&1
\end{array}
\right]v.
$$
The $3\times 3$ matrix above defines a metric at $p$. Varying $p$ we obtain a Riemannian metric $\dot{\cdot}{\cdot}$, since each matrix entry is differentiable. The vectors $(1,0,0)$, $(0,1,p_x)$,  and $(0,0,1)$ form an orthogonal basis at $p$. Also, the volume form of $Nil$ coincides with the standard one from $\R^3$, since the metric determinant is equal to one.  

The geodesic flow on $Nil$ admits a solution~\cite{szilagyi2003curvature}. A ray $\gamma(t)=(x(t),y(t),z(t))$ starting at $(0,0,0)$ in the tangent direction $v=(c\cos(\alpha),c\sin(\alpha), w)$ is given by:
$$
\begin{array}{l}
\begin{array}{l}
\displaystyle x(t)=\frac{c}{w}(\sin(wt+\alpha)-\sin(\alpha)),
\end{array}
\vspace{0.1cm}
\\
\begin{array}{l}
     \displaystyle
y(t)=-\frac{c}{w}(\cos(wt+\alpha)-\cos(\alpha)),
\end{array}
\vspace{0.cm}
\\
\begin{array}{ll}
\displaystyle
z(t)= t(w +\frac{c^2}{2w})&\hspace{-0.3cm}\displaystyle-\frac{c^2}{4w^2}(\sin(2wt+2\alpha)-\sin(2\alpha))
\vspace{0.2cm}
\\
&\hspace{-0.3cm}\displaystyle+\frac{c^2}{2w^2}(\sin(wt+2\alpha)-\sin(2\alpha)-\sin(tw)).
\end{array}
\end{array}
$$
\newpage

To compute a geodesic $\beta(t)$ starting at $p$ in the direction $v$, we first translate the initial conditions to the origin and then compute the geodesic using the solution above. We translate this back to the desired position.
Note that geodesics joining two points in $Nil$ are in general not unique.

For a compact manifold $M$ modeled by $Nil$ consider a discrete group $\Gamma$ generated by the ``translations'' in the direction of axis $x$, $y$, and $z$: $\Phi_1(p)=(x+1, y, y+z)$, $\Phi_2(p)=(x,y+1,z)$, and $\Phi_3(p)=(x,y,z+1)$. The manifold $M$ inherits the geometry of $Nil$. For each fixed $x$ we obtain a two-dimensional torus and $M$ is foliated by tori. The unit cube is the fundamental domain. 
We set the scene inside the cube and ray trace it using the concepts from Section~\ref{sec:shader}. Each time a ray intersects a cube face we update it using $\Gamma$. 

\begin{figure}[ht]
	\centering
	\includegraphics[width=0.9\columnwidth]{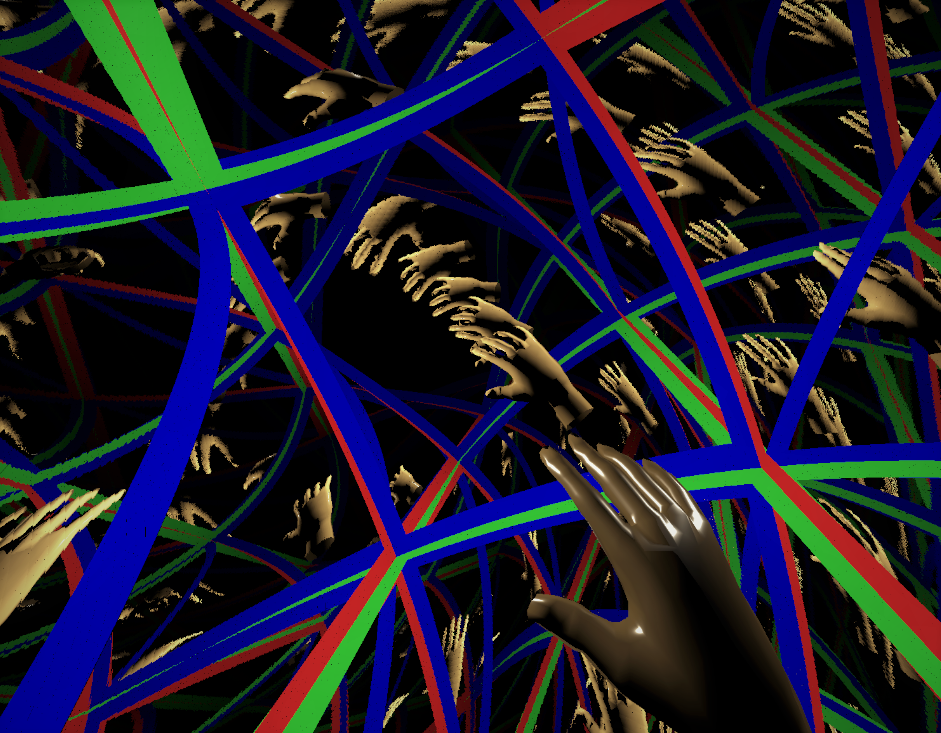}
	\caption{Inside view of a $Nil$ manifold. A hand and the cube's edges compose the scene. The face pairing makes the rays iterate giving a tessellation of $Nil$ by cubes.}
	\label{fig:nil}
\end{figure}

Figure~\ref{fig:nil} gives an inside view of $M=Nil/\Gamma$. The scene is composed of thickening in 2D of the boundary of the fundamental domain faces.
Opposite faces get the same shading: red, green, and blue for the faces parallel to the axis, $x$, $y$, and $z$, respectively.
In the images, the lines shaded with red and green are extended to infinity as straight lines. These lines are geodesics of $Nil$ and are perpendicular to the plane $yz$; lines with such property foliate the space $Nil$.

The three-dimensional Heisenberg group $H$ provides a basic example of \emph{contact geometry} (see \cite[Sections 0.2 and 0.3]{Gromov96}). Here the idea is to define a distance on $H$ restricted to a polarization given by the left translates of the $(x,y)$-plane. The distance between two points is the infimum of the length of the paths whose tangents belong to the polarization. Because the $Nil$ space is curved in a particular way, it appears that any two points can be joined by such a path. The resulting metric is not Riemannian, it is an example of so called Carnot--Carath\'{e}odory metric and it is very different from the Riemannian metric of $Nil$. This construction gives the simplest non-trivial example of a contact Carnot--Carath\'{e}odory metric space which is arguably the most important at the same time. We refer to a beautiful monograph \emph{``Carnot--{C}arath\'{e}odory spaces seen from within''} by M.~Gromov for an introduction to this subject \cite{Gromov96}. Inside view of the  Carnot--Carath\'{e}odory spaces is an interesting project for the future investigation.


\subsubsection{$\USL$ space}\label{sec:sl2_space}

The $\USL$ geometry is similar to $Nil$, but it is now an $\R$-bundle over $\H^2$. This
geometry is constructed from the Lie group $\SL$.
In \cite{nilsolsl2} the space $\USL$ is described in more details focusing on the technical aspects of visualization.

We follow the notation of Gilmore~\cite{gilmore2008lie}. The \textit{special linear group} $\SL$ consisting of all $2\times 2$ matrices with unit determinant is a Lie group: the product of two matrices with unit determinant has unit determinant, the same for the inverse matrix.

To understand the hyperbolic nature of $\SL$ observe that the elements of $\SL$ are matrices $\begin{psmallmatrix}a & b\\c & d\end{psmallmatrix}$
such that $ad-bc=1$. Hence $\SL$ is a $3$-manifold whose embedding in $\R^4$ is given by $\{(a,b,c,d)\in \R^4|\,\, ad-bc=1\}$. Rewriting $ad-bc=1$, we get:
$$\left(\frac{a+d}{2}\right)^2-\left(\frac{a-d}{2}\right)^2+\left(\frac{b-c}{2}\right)^2-\left(\frac{b+c}{2}\right)^2=1,$$
which is an equation of a three-dimensional hyperbola in $\R^4$.

There is also an identification of $\SL$ with $\H^2\times \S^1$; see~\cite{nilsolsl2} for more details. That is, $\SL$ is not simply connected, which implies that it is not a model geometry. However, the \textit{universal cover} $\USL$ of $\SL$ is a model geometry (see~\cite{thurston1979}).
We focus on the visualization of $\SL$ since the geometries are locally isometric.

We use the parameterization of a neighborhood of identity in $\SL$ from \cite{gilmore2008lie}:
\begin{equation}\label{eq:sl2_parameterization}
	\textbf{x}(x,y,z)=\left[
	\begin{array}{cc}
		1+x&y
		\vspace{.3cm}\\
		z&\displaystyle\frac{1+yz}{1+x}
	\end{array}
	\right].
\end{equation}
Observe that $\textbf{x}(0,0,0)$ is the identity of $\SL$, and that in the plane $x=-1$ the map $\textbf{x}$ is not defined. We use $\textbf{x}$ to push-back the metric of $\SL$ to $\R^3$. 

We now construct a metric in the $\SL$.
The element $e = \begin{psmallmatrix}1 & 0\\0 & 1\end{psmallmatrix}$ is the identity of $\SL$. Let $T_e\SL$ be the tangent space at $e$ with the well-known scalar product
$\dot{u}{v}_e=\mbox{Trace}(u\cdot v)$
between two tangent vectors $u$ and $v$~\cite{gilmore2008lie}.
 As in $Nil$ geometry,
we extend it to a Riemannian metric using left multiplication.

Using the above Riemannian metric we obtain the geodesic flow below, which can be solved using Euler's method (see~\cite{nilsolsl2}): 

\newpage

\begin{equation}\label{eq:geodesic_equation_bundle_SL2}
	\left\{
	\begin{array}{ll}
		x'_k& = y_k, \,\, k=1,2,3,   \\[0.0cm]
		y'_1& = \displaystyle\frac{(1+p_yp_z)y_1^2}{1+p_x}-p_zy_1y_2-p_yy_1y_3+(1+p_x)y_2y_3,
		\\[0.25cm]
		y'_2& = \displaystyle\frac{(1+p_yp_z)p_yy_1^2}{(1+p_x)^2}-\frac{p_zp_y}{1+p_x}y_1y_2-\frac{p_y^2}{1+p_x}y_1y_3+p_yy_2y_3,
		\\[0.25cm]
		y'_3& = \displaystyle\frac{(1+p_yp_z)p_zy_1^2}{(1+p_x)^2}-\frac{p_z^2}{1+p_x}y_1y_2-\frac{p_yp_z}{1+p_x}y_1y_3+p_zy_2y_3.
	\end{array}
	\right.
\end{equation}

In \cite{divjak2009geodesics, marenitch2008geodesic}, the authors provide compact formulas for the geodesics of $\SL$, however, they use a different parameterization of the space. It could be interesting to explore those formulas to increase performance in the visualizations.

\begin{figure}[ht]
	\begin{center}
	\includegraphics[width=0.9\columnwidth]{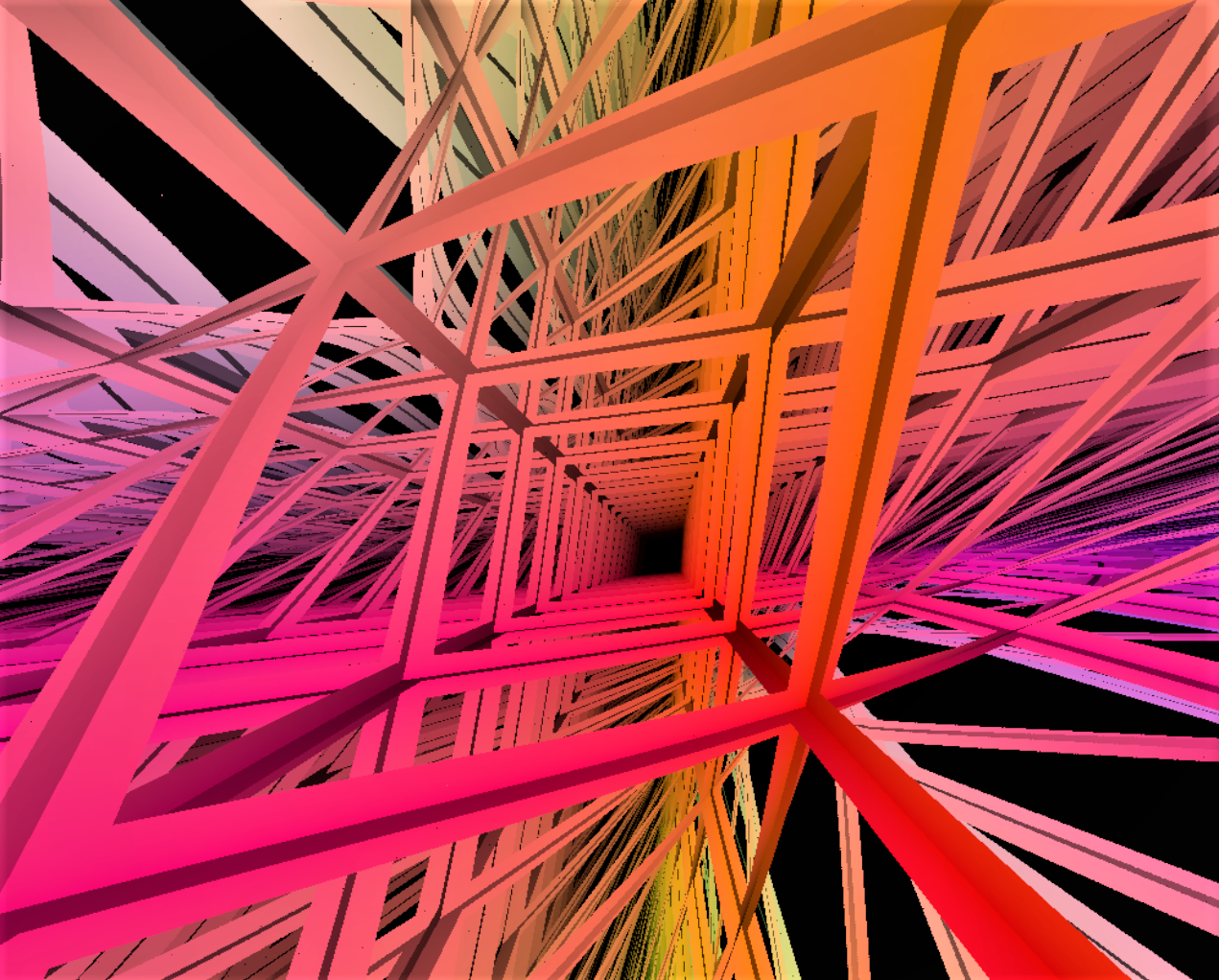}
	\caption{Inside view of $\SL$. The scene is a grid in $\R^3$ deformed by the $\SL$~metric.}
	\label{fig:sl2}
	\end{center}
\end{figure}

Figure~\ref{fig:sl2} presents a visualization of $\SL$.
The image is a visualization of an quotient manifold of $\USL$ by an infinite cyclic discrete group. It is a non-compact and infinite volume manifold. The picture presents an inside view of the special linear group $\SL$ which has the geometry modeled by $\USL$.
The image refers to an immersive visualization of a grid defined in the domain of the parameterization $\textbf{x}(x,y,z)$ of $\SL$. The parameterization distorts the grid following the geometry of $\SL$. The choice of the RGB colors was done empirically for the $(x,y,z)$ coordinates of the hit points.

\subsubsection{$Sol$ space}

$Sol$ is the least symmetric among the eight Thurston's geometries. It is a plane
bundle over the real line.
Its geometry comes from a solvable Lie group $Sol$. For more details see~\cite[Section 12.7]{martelli2016introduction}.

The \textit{$Sol$ space} ($Sol$) is an example of a Lie group which consists of all matrices 
$$\small\left[
\begin{array}{lll}
e^{z}&0&x\\
0&e^{-z}&y\\
0&0&1
\end{array}
\right]$$
with the multiplication operation. Clearly, $Sol$ is diffeomorphic to $\R^3$.

Let $(x,y,z)$ and $(x',y',z')$ be two elements in $Sol$. Their multiplication has the form:
$$
(x,y,z)\cdot (x',y',z')=(x'e^{z}+x,y'e^{-z}+y,z+z'),
$$
which is the sum of the element coordinates controlled by an additional term in the first coordinates. To endow $Sol$ with a geometry we consider the Euclidean metric in the tangent space at the origin and extend it by left multiplication. 
After some computations we get the scalar product of two tangent vectors $u$ and $v$ at $p$:
$$\dot{u}{v}_p=u^T\left[
\begin{array}{ccc}
e^{2p_z}&0&0\\
0&e^{-2p_z}&0\\
0&0&1
\end{array}
\right]v.$$
The matrix above defines a metric at $p$. Varying $p$ we obtain a Riemannian metric $\dot{\cdot}{\cdot}$, since each matrix entry is differentiable. The volume form of $Sol$ coincides with the standard one from $\R^3$, since the determinant of the above matrix is one. 

Using the above metric we obtain the \textit{geodesic flow} of the $Sol$ geometry:
\begin{equation}\label{eq:geodesic_equation_bundle_sol}
    \left\{
    \begin{array}{ll}
         x'_k& = y_k, \,\, k=1,2,3,   \\
          y'_1& = -2y_1y_3,
          \\
          y'_2& = 2y_2y_3,
          \\
          y'_3& = e^{2p_z}y_1^2-e^{-2p_z}y_2^2.
    \end{array}
    \right.
\end{equation}
Unfortunately, there is no solution for this system in terms of elementary functions~\cite{szilagyi2007Sol}. 
Troyanov~\cite{troyanov1998horizon} obtained a formula for geodesics in $Sol$, however, it contains many coefficients that can not be computed in a closed formula. He classified $Sol$ geodesics in classes of equivalence, the \textit{horizons} of $Sol$.

Similarly to the case of $Nil$, geodesics joining two points in $Sol$ geometry are not~unique.

\newpage

Let $\Gamma$ be a discrete group generated by the translations $\Phi_1$ and $\Phi_2$ along $x$- and $y$-axes, and a transformation $\Phi_3$: 
\begin{align*}
\Phi_1(p) &=(x+1, y, z),\\ 
\Phi_2(p) &=(x,y+1,z),\\
\Phi_3(p) &=(x\cdot e^{-2\ln\phi},y\cdot e^{2\ln\phi},z+2\ln\phi);
\end{align*}
where $p=(x,y,z)$ and $\phi$ is the golden ratio number. 

The manifold $M=Sol/\Gamma$ inherits the geometry of $Sol$, and for each fixed $z_0$ it contains a $2$-torus which is the quotient of $\R^2\times \{z_0\}$ by the discrete group generated by $\Phi_1$ and $\Phi_2$. Thus $M$ is foliated by the $2$-tori. The parallelepiped $D\times [0,2\ln\phi)$ is a fundamental domain of $\Gamma$, where $D$ is the unit square.

\begin{figure}[ht]
	\centering
	\includegraphics[width=0.9\columnwidth]{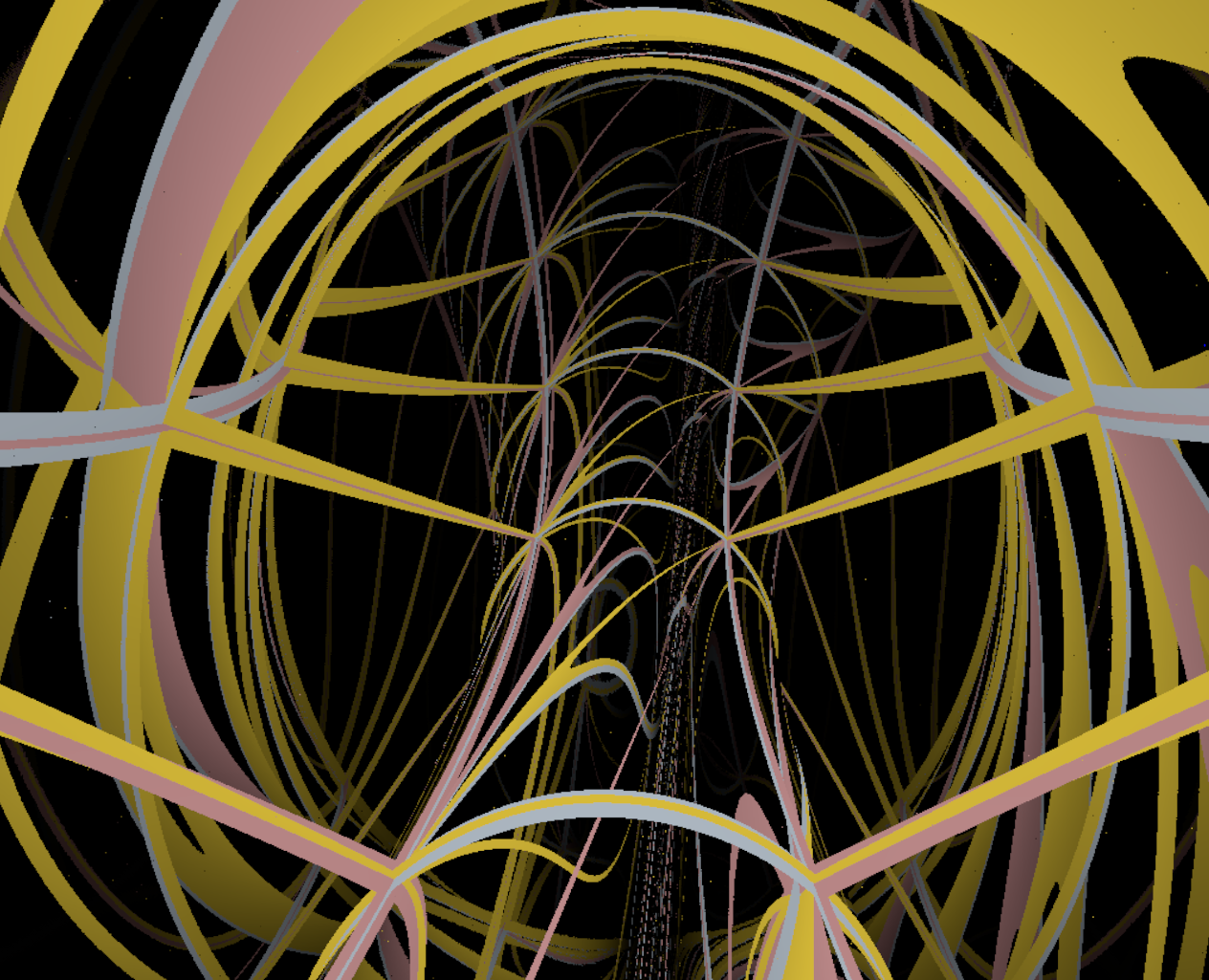}
	\caption{Inside view of a $Sol$ manifold. The scene is composed by the cube's edges. The face pairing makes the rays that leave a face return from its opposite, tessellating $Sol$.}
	\label{fig:sol}
\end{figure}

The scene is set inside the parallelepiped and ray traced using $Sol$ rays. Each time a ray hits a parallelepiped face it is updated by the discrete group $\Gamma$.
As in the $Nil$ section, a Riemannian shading is defined for $Sol$. The scene was created by thickening the boundary of the fundamental domain faces and the opposite faces get the same shading.
Figure~\ref{fig:sol} presents the visualization of the manifold $M=Sol/\Gamma$. 

\section{Final remarks}\label{sec-final_remarks}

This expository paper provides a brief overview of the history, main definitions, and important results regarding Thurston's geometries. We focused on an intuitive presentation of the geometrization conjecture which, roughly speaking, states that any $3$-manifold can be cut into nice geometric pieces, each modeled by one of the eight  Thurston's geometries. This fundamental conjecture, which is now a theorem due to Perelman, motivated computer graphics algorithms for immersive visualization of such geometric structures. In the paper we present the images that were rendered by using the recent Riemannian ray tracing algorithm.

\begin{figure}[!ht]
    \centering
    \includegraphics[width=0.95\linewidth]{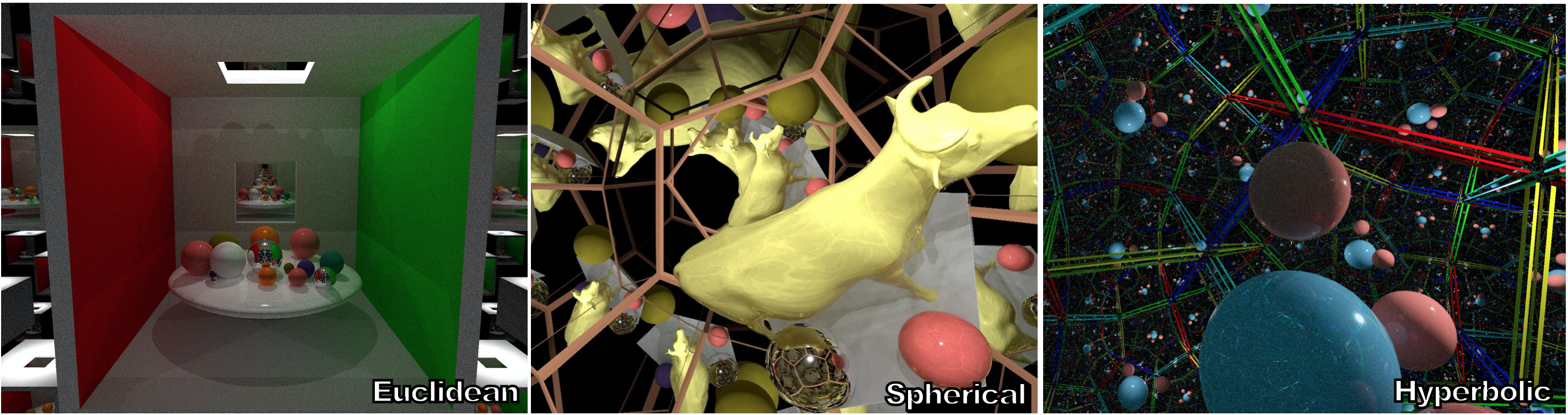}
    \vspace{-0.3cm}
    \caption{``Photorealistic'' inside views of the flat torus, Poincaré sphere and hyperbolic mirrored dodecahedron. For more details on these visualizations and the techniques used for rendering the images see~\cite{global_illumination}.}
    \label{fig:global_illumination}
\end{figure}

We believe that the development of Computer Graphics techniques in Riemannian geometry, such as the Riemannian ray tracing, could be an ally in mathematical research in low-dimension geometry and topology. 
Being inside these abstract landscapes by means of the Virtual Reality devices could inspire new ideas or help us to gain more intuition on these geometries.
Several ongoing projects have been using the Virtual Reality to interactively visualize the eight Thurston geometries~\cite{nilsolsl2,velho2020immersive, 33cmb,weeks2020billiards,weeks2021body, coulon2020ray}.
Furthermore, very recently, a \textit{global illumination} algorithm was introduced to produce ``photorealistic'' inside views of Thurston's classical geometries (see Figure~\ref{fig:global_illumination}). We refer to  \cite{global_illumination} for more details on this project.

\begin{figure}[!ht]
    \centering
        \includegraphics[width=0.225\columnwidth]{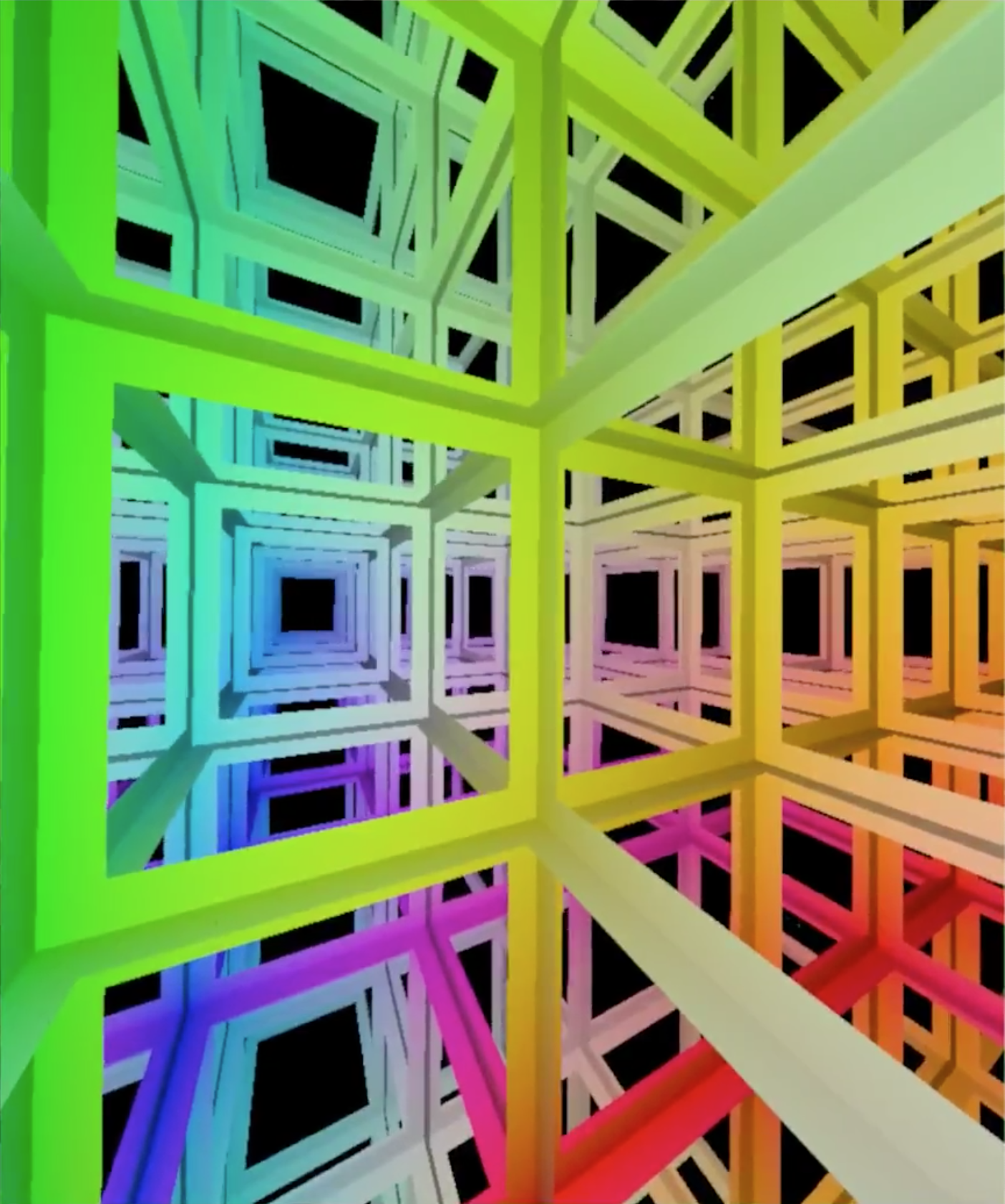}
        \includegraphics[width=0.225\columnwidth]{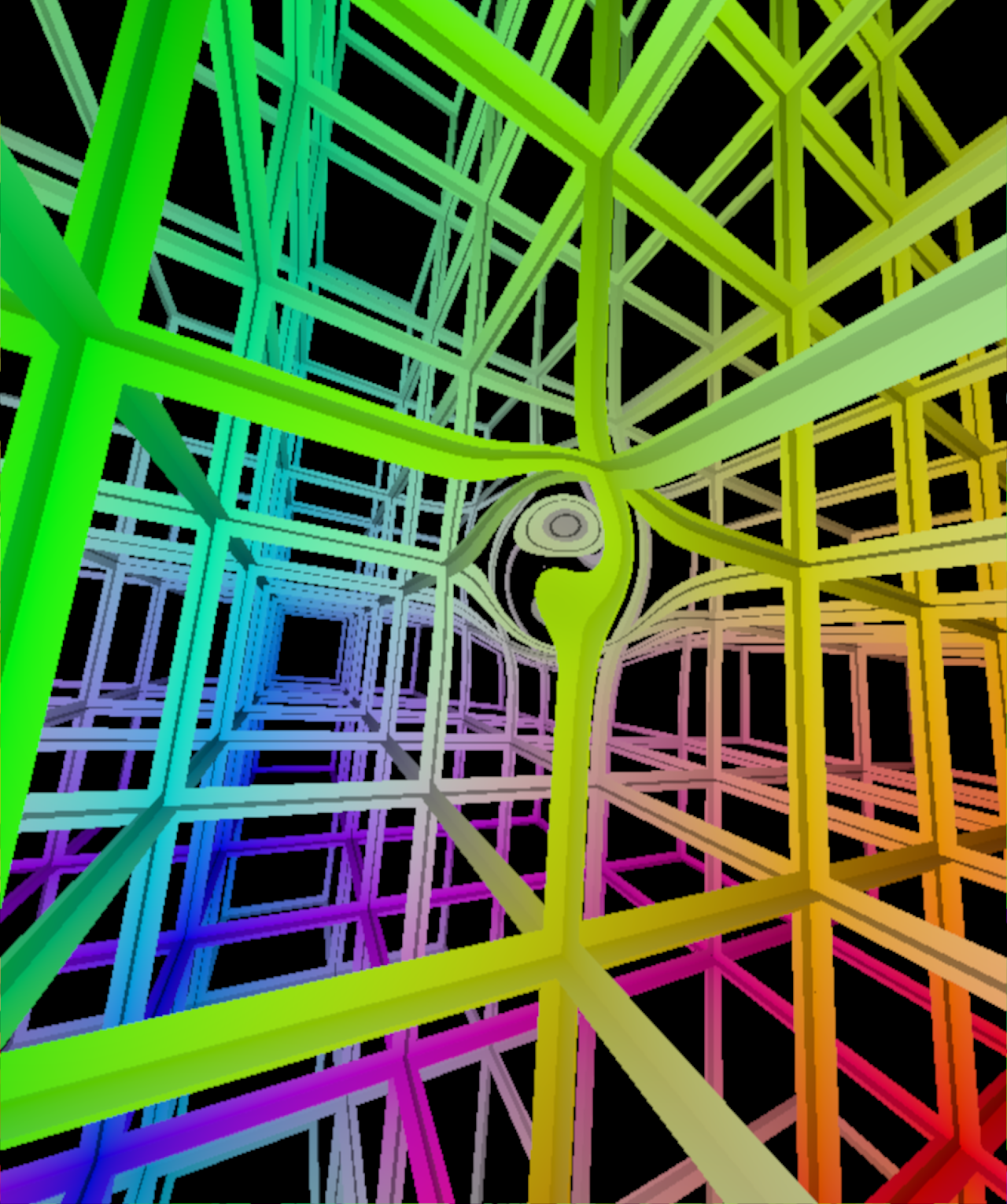}
        \includegraphics[width=0.225\columnwidth]{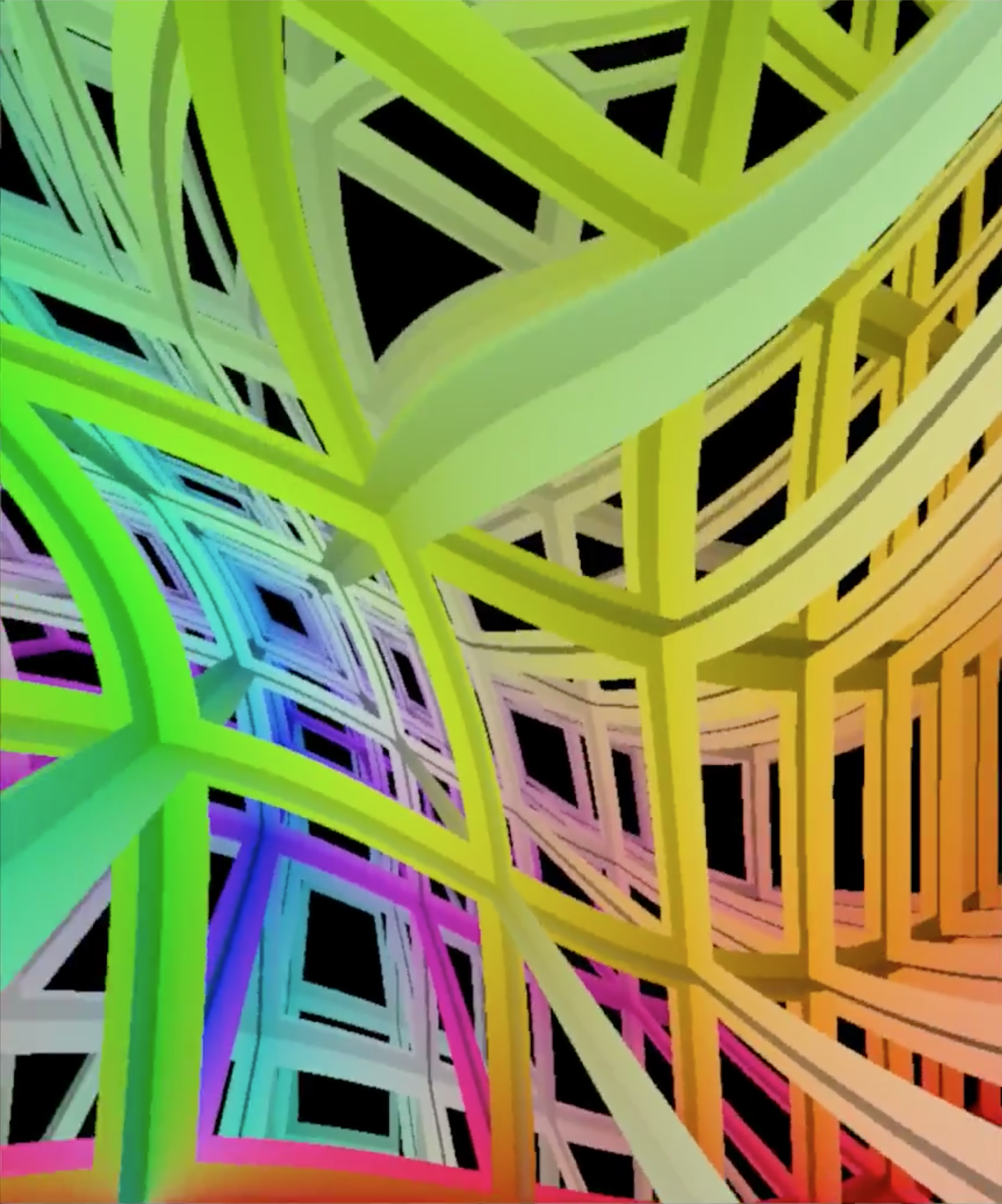}
        \includegraphics[width=0.225\columnwidth]{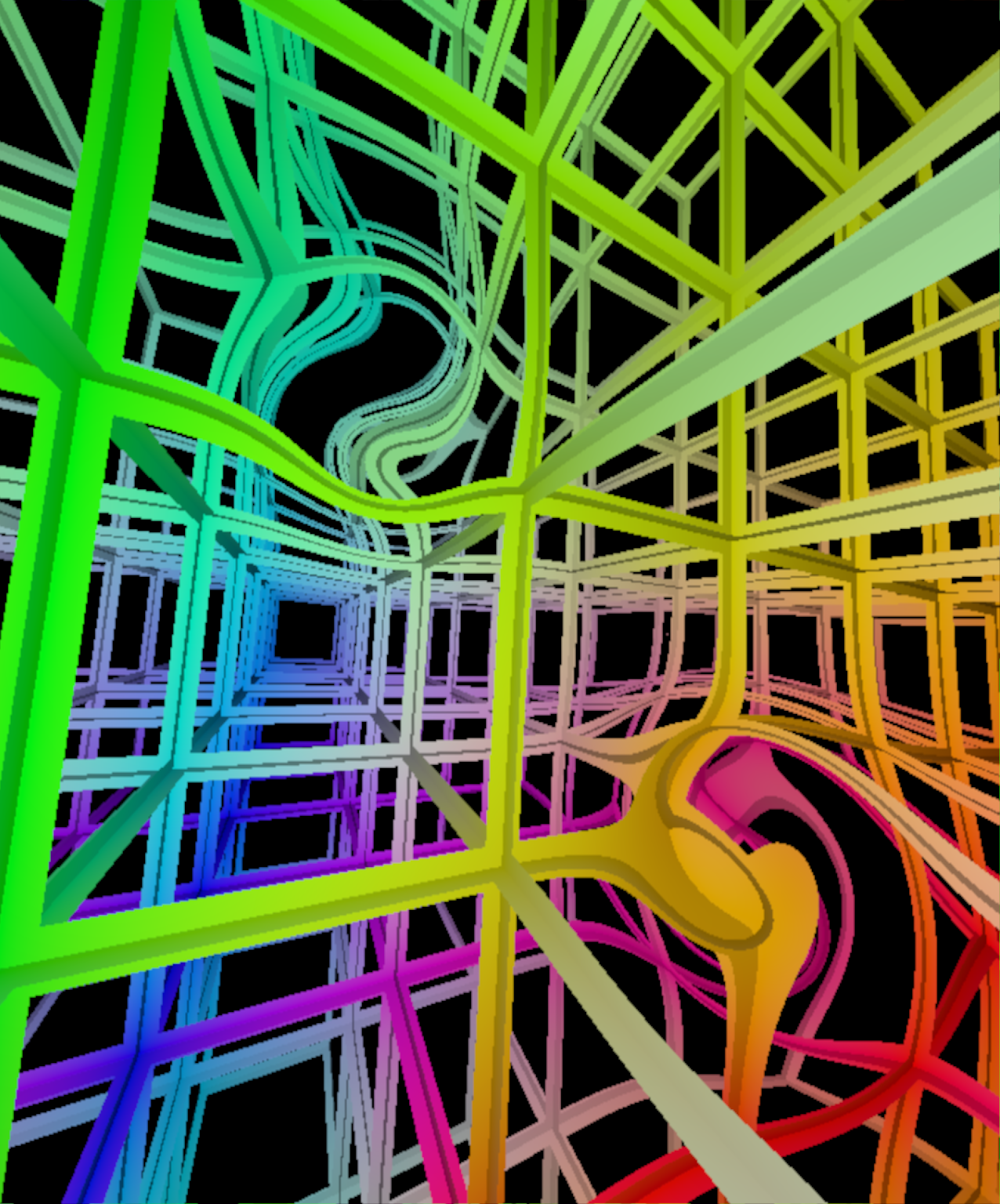}
    \vspace{-0.3cm}
    \caption{The figures provide immersive visualizations of three-dimensional graphs of Gaussians. We set our scenes in $\R^3$ and we push-back the metric of the graph to $\R^3$ to ray trace a regular grid. From the left to right: no deformation, local deformation, global deformation, two local deformations using the sum of Gaussians.}
    \label{fig:deformations}
\end{figure}

On the other hand, Riemannian geometry tools could help us in development of new computer graphics algorithms.
For example, a path for future work would be focusing on time/user-dependent Riemannian metric constructions on the ambient space to explore the special effects like warping, mirages~\cite{stam1996ray}, and scene deformation~\cite{barr1986ray}. 
In~\cite{novello2020design}, the authors used graphs of functions to design Riemannian metrics on $\R^3$, allowing modeling of expressive effects (see Figure~\ref{fig:deformations}).
We expect that curved rays can advance the state of the art in many areas, not restricted to rendering only. 

In conclusion, we can mention some other interesting projects that follow the non-Euclidean visualization program presented in this work.
The visualization of surfaces embedded in the eight Thurston geometries is a feasible problem. 
For example, being inside the $3$-sphere would allow us to visualize the flat torus, which is not possible in the Euclidean $3$-space. More general visualization of minimal and constant curvature surfaces is a natural continuation of this proposal, we expect that in this context the product geometries of Thurston will play a prominent role.
Inside view of the other three-dimensional metric spaces, such as the Carnot--Carath\'{e}odory spaces briefly discussed in Section~\ref{sec:Nil}, is another interesting project for the future investigation.
\bibliographystyle{amsplain}
\bibliography{ray3d}

\end{document}